\documentclass{article}

\usepackage{graphicx}
\usepackage{url}
 \usepackage[pdftex]{color}
\usepackage{amsmath}
\usepackage{amsfonts} 
\usepackage{array}
\usepackage{cancel}
\usepackage[utf8]{inputenc}
\usepackage{amsmath}
\usepackage{amssymb}
\usepackage{makeidx}
\usepackage{multicol}
\usepackage[bottom]{footmisc}
\usepackage{graphicx}
\usepackage{color}
\usepackage{epstopdf}
\usepackage{algorithm}
\usepackage{algcompatible}
\usepackage{calc}

\newtheorem{teor}{Theorem}[section]
\newtheorem{cor}{Corollary}[section]
\newtheorem{propo}{Proposition}[section]
\newtheorem{lemm}{Lemma}[section]

\newtheorem{probl}{Problem}
\newtheorem{rem}{Remark}
\newtheorem{approxi}{Approximation}
\newtheorem{assump}{Assumption}

\newcommand{\app}{\begin{approxi}}
\newcommand{\eapp}{\end{approxi}}
\newcommand{\ass}{\begin{assump}}
\newcommand{\eass}{\end{assump}}
\newcommand{\teo}{\begin{teor}}
\newcommand{\eteo}{\end{teor}}
\newcommand{\corr}{\begin{cor}}
\newcommand{\ecorr}{\end{cor}}
\newcommand{\pro}{\begin{propo}}
\newcommand{\epro}{\end{propo}}
\newcommand{\lemma}{\begin{lemm}}
\newcommand{\elemma}{\end{lemm}}
\newcommand{\pb}{\begin{probl}}
\newcommand{\epb}{\end{probl}}
\newcommand{\df}{\begin{defn}}
\newcommand{\edf}{\end{defn}}
\newcommand{\rema}{\begin{rem}}
\newcommand{\erema}{\end{rem}}
\newcommand{\qed}{\hfill$\blacksquare$}

\newcommand{\al}[1]{\begin{align} #1 \end{align}}
\newcommand{\nn}{\nonumber}

\renewcommand{\d}{\mathrm d}
\newcommand{\Ac}{ \mathcal{A}}

\newcommand{\Dc}{ \mathcal{D}}

\newcommand{\Fc}{ \mathcal{F}}
\newcommand{\Gc}{ \mathcal{G}}

\newcommand{\Ic}{ \mathcal{I}}

\newcommand{\Kc}{ \mathcal{K}}

\newcommand{\Nc}{ \mathcal{N}}
\newcommand{\Oc}{ \mathcal{O}}

\newcommand{\Sc}{ \mathcal{S}}

\newcommand{\Uc}{ \mathcal{U}}
\newcommand{\Vc}{ \mathcal{V}}
\newcommand{\Wc}{ \mathcal{W}}
\newcommand{\Xc}{ \mathcal{X}}
\newcommand{\Yc}{ \mathcal{Y}}

\newcommand{\Ns}{ \mathbb{N}}

\newcommand{\Rs}{ \mathbb{R}}

\newcommand{\diag}{ \operatorname{diag}}
\newcommand{\tpl}[1]{ \operatorname{Tpl}(#1)}
\date{}
\begin{document}

\title{A second-order generalization\\ of TC and DC kernels}

\author{Mattia~Zorzi \thanks{M. Zorzi are with   the Department of Information
Engineering, University of Padova, Padova, Italy; e-mail: {\tt\small zorzimat@dei.unipd.it} (M. Zorzi).}}

\maketitle

\begin{abstract}Kernel-based methods have been successfully introduced in system identification to estimate the impulse response of a linear system.  Adopting the Bayesian viewpoint, the impulse response is modeled as a zero mean Gaussian process whose covariance function (kernel) is estimated from the data. The most popular kernels used in system identification are the tuned-correlated (TC), the diagonal-correlated (DC) and the stable spline (SS) kernel. TC and DC kernels admit a closed form factorization of the inverse. The SS kernel induces more smoothness than TC and DC  on the estimated impulse response, however, the aforementioned property does not hold in this case. In this paper we propose a second-order extension of the TC and DC kernel which induces more smoothness than TC and DC, respectively, on the impulse response and a generalized-correlated kernel which incorporates the TC and DC kernels and their second order extensions. Moreover, these generalizations  admit a closed form factorization of the inverse and thus they allow to design efficient algorithms for the search of the optimal kernel hyperparameters. We also show how to use this idea to develop higher oder extensions. Interestingly, these new kernels belong to the family of the so called exponentially convex  local stationary kernels: such a property allows to  immediately analyze the frequency properties induced on the estimated impulse response by these kernels.    \end{abstract}

\section{Introduction}\label{sec_intro}
Linear system identification problems are traditionally addressed 
by using Prediction Error Methods (PEM), see \cite{LJUNG_SYS_ID_1999,SODERSTROM_STOICA_1988}.
Here, the best model is chosen over a fixed parametric model class (e.g. ARMAX, OE, Box-Jenkins). This approach, however, has two issues: first, the parametrization of the predictor is nonlinear which implies that the minimization of the squared prediction error leads to a non-convex optimization problem; second, we have to face a model selection problem (i.e. order selection) which is usually performed by AIC and BIC criteria 
\cite{AKAIKE_1974,SCHWARZ_1978}.

Regularized kernel-based methods have been recently proposed in system identification in order to overcome the aforementioned limitations, see  \cite{PILLONETTO_DENICOLAO2010,EST_TF_REVISITED_2012,KERNEL_METHODS_2014}. Here, we search the candidate model, described via the predictor impulse response, in an infinite dimensional nonparametric model class with the help of a penalty term. Adopting the Bayesian viewpoint, this is a Gaussian process regression problem \cite{RASMUSSEN_WILLIAMNS_2006}: the impulse response is modeled as a Gaussian process with zero mean and with a suitable covariance function, also called kernel \cite{wahba1990spline}. The latter encodes the a priori knowledge about the predictor impulse response. For instance, the impulse response should be Bounded Input Bounded Output (BIBO) stable and with a certain degree of smoothness.

The most popular kernels are the tuned-correlated (TC), the diagonal-correlated (DC), and the stable-spline (SS), see \cite{EST_TF_REVISITED_2012,PILLONETTO_DENICOLAO2010}. All these kernels encode the BIBO stability property.
Regarding the smoothness, SS is the one inducing more smoothness on the impulse response. It is worth noting that many other extensions can be obtained, see for  instance  \cite{CHEN2018109,dinuzzo2015kernels,ZORZI2018125,zorzi2020new}. All these kernels depend on few hyperparameters that are learnt from the data by minimizing the so called negative log-marginal likelihood. This task is computationally expensive especially in the case we want to estimate high dimensional models, e.g. the case of dynamic networks, see  \cite{BSL,CHIUSO_PILLONETTO_SPARSE_2012,BSL_CDC}.

To reduce the computational complexity different strategies have been proposed, see  \cite{carli2012efficient,chen2018regularized, chen2021semiseparable,chen2013implementation}. In particular, if the  kernel matrix admits a closed form expression for Cholesky factor of its inverse matrix (ant thus also its determinant), then the evaluation of the marginal likelihood can be done  efficiently  \cite{chen2013implementation}. While it is possible to derive these closed form expressions for TC and DC,  see \cite{CARLI2014,7495008}, this is not possible for SS. It is worth noting that an efficient algorithm for the SS kernel has been proposed in  \cite{chen2021semiseparable}. The latter, however, can be used only in the case that the input of the system has a prescribed structure, e.g. it cannot be used in the case we collect the data of a system which is in a feedback configuration.

The aim of this paper is to introduce a second-order generalization of the TC and DC kernel exploiting the filter-based approach proposed in \cite{marconato2017filter}. These extensions induce more smoothness than TC and DC, respectively. We also introduce a generalized-correlated kernel which incorporates the DC, TC kernels and their second order extensions. Moreover, we show that  they admit a closed form expression for the Cholesky of its inverse matrix. Thus, these kernels allow to design an efficient algorithm for the search of the optimal hyperparameters. It is worth noting that SS is the second-order extension of the TC kernel derived in the continuous time. In contrast, the extension that we propose here is derived in the discrete time.
 Numerical experiments showed that the new second-oder TC kernel represents an attractive  alternative to SS because it leads to an estimation algorithm which outperforms the one using SS (even in the case that the computation of the Cholesky factorization of the kernel exploits the fact that SS is extended 2-semiseparable) in terms of computational complexity, while the second-order TC and SS are similar in terms of estimation performance. This idea can be also used to higher order extensions and also to generalize the high frequency kernel proposed in \cite{6160606}. Interestingly, all these new kernels are exponentially convex local stationary (ECLS),  \cite{CHEN2018109,ZORZI2018125}. Such a property allows to easily understand the frequency properties of their stationary parts.

The outline of the paper is as follows. In Section \ref{sec_PEM} we briefly review the kernel-based PEM method as well as the TC, DC and SS kernels. Section 
\ref{sec_TC2} introduces the second-order extension for the TC kernel, while Section \ref{sec_DC2} the one for the DC kernel. In Section \ref{sec_GC} we introduce the generalized-correlation kernel.
In Section \ref{sec_ME} we derive the closed form expressions for these kernels. In Section \ref{sec_high} we extend this idea to higher order generalizations.  In Section \ref{sec_freq} we show that these kernels are ECLS and we analyze the stationary part of these kernels in the frequency domain. Finally, we draw the conclusions in Section  \ref{sec_conc}.

{\em Notation}. $\Sc_T$, with $T\leq \infty$, denotes the cone of positive definite symmetric matrices of dimension $T\times T$. Infinite dimensional matrices, i.e. matrices having an infinite number of columns and/or rows, are denoted using the calligraphic font, e.g. $\Kc$, while finite dimensional ones are denoted using the normal font, e.g. $K$. Given $\Fc\in\Rs^{p\times \infty}$ and $\Gc\in\Rs^{\infty\times m}$, the product  $\Fc\Gc$ is understood as a $p\times m$ matrix whose entries are limits of infinite sequences \cite{INFINITE_MATRICES}. Given $K\in\Sc_T$, $[K]_{t,s}$ denotes the entry of $K $in position $(t,s)$, while $[K]_{:,t}$ and $[K]_{t,:}$ denotes the $t$-th column and row, respectively, of $K$. Given $K\in\Sc_T$, $\|v\|_{K^{-1}}$ denotes the weighted Euclidean norm of $v$ with weight $K^{-1}$.
Given $v\in\Rs^T$, $\tpl v$ denotes the lower triangular $T\times T$ Toeplitz matrix whose first column is given by $v$, while $\diag(v)$ denotes the diagonal matrix whose main diagonal is $v$.

\section{Kernel-based PEM method} \label{sec_PEM}
Consider the  model 
\al{\label{mod}y(t)=\sum_{k=1}^\infty g(k) u(t-k)+e(t), \quad t=1\ldots N}
where $y(t)$,  $u(t)$, $g(t)$ and $e(t)$ denote the output, the input, the impulse response of the model and a zero-mean white Gaussian noise with variance $\sigma^2$, respectively.  We can rewrite model (\ref{mod}) as
\al{y=\Ac g+e\nn}
where $y=[\,y(1)\ldots y(N)\,]^\top\in\Rs^N$, $e$ is defined likewise, $\Ac^{N\times \infty}$ is the regression matrix whose entries are defined by $u(t)$ with $t=1\ldots N$, $g=[\,g(1) \; g(2)\ldots \,]^\top\in\Rs^\infty$. We want to estimate the impulse response $g$ given the measurements $\{y(t),u(t)\}_{t=1}^N$. Such a problem is ill-posed because we have a finite number of measurements while $g$ contains infinite parameters. The latter can be made well-posed assuming that $g\sim \Nc(0,\lambda \Kc(\eta))$ where $\Kc(\eta)\in\Sc_\infty$ is the kernel function and $\eta$ is the vector of hyperparameters characterizing the kernel; in this  
way, the minimum variance estimator of $g$ is:
\al{\label{def_RELS}\hat g=\underset{g\in\Rs^\infty}{\mathrm{argmin}} \|y-\Ac g \|^2+\frac{\sigma^2}{\lambda} \| g\|^2_{\Kc(\eta)^{-1}}}
where $\lambda>0$ denotes the regularization parameter. It is worth noting that the above problem admits a closed form solution.
Moreover, $\Kc(\eta)$ encodes the a priori information that we have on the impulse response. 

The aforementioned problem can be formulated as a finite dimensional problem. Indeed, $g$ can be truncated, obtaining a finite impulse response of length $T$; the corresponding kernel matrix $K(\eta)\in\Sc_T$ is defined as $[K(\eta)]_{t,s}=[\Kc(\eta)]_{t,s}$ for $t,s=1\ldots T$ and the regression matrix $A\in\Rs^{N\times T}$ is given by the first $T$ columns of $\Ac$.
Such a truncation, with $T$ sufficiently large, does not introduce a bias, because $g$ decays to zero. The so called hyperparameters $\lambda$ and $\eta$ are estimated by minimizing numerically the negative log-marginal likelihood 
\al{\label{marginal}\ell(y; \lambda, \eta):=\log &\det(\lambda A K(\eta)A^\top+\sigma^2 I)\nn\\ & +y^\top( \lambda A K(\eta)A^\top+\sigma^2 I)^{-1}y.} 
In what follows, we will drop the dependence on $\eta$ for kernels in order to ease the notation. 

\subsection{Diagonal and correlated kernels: an overview}
We briefly review the most popular kernels used in system identification, see \cite{EST_TF_REVISITED_2012} for a more complete overview. The simplest kernel is diagonal and encodes the a priori information that $g$ should decay to zero exponentially:
\al{\label{kernelDI} \Kc_{DI} =  \diag(\beta,\beta^2,\ldots, \beta^t,\ldots)}
where $\eta=\beta$ and $0<\beta<1$. Indeed, the penalty term $\|g\|^2_{\Kc_{DI}^{-1}}$ is the squared norm of the weighted impulse response \al{h=[\, h_1 \; h_2 \ldots h_t\ldots \,]^\top, \quad h_t=\beta^{-t/2}g_t\nn} which amplifies in an exponential way the coefficients $g_t$ as $t$ increases. The tuned-correlated (TC, also called first-order stable spline) kernel embeds also the a priori information that $g$ is smooth: 
\al{\label{kernelTC} [\Kc_{TC}]_{t,s}  =  \beta^{\max(t,s)}}
where $\eta=\beta$ and $0<\beta<1$. The smoothness property can be justified as follows. It is well known that  
\al{ \Kc_{TC}  =(1-\beta)  (\Fc\Dc\Fc^T)^{-1}\nn}
where 
\al{\Fc&=\tpl{1,-1,0,\ldots }\nn\\
\Dc&=\diag(\beta^{-1},\beta^{-2},\ldots, \beta^{-t},\ldots);\nn} then, $\Fc^\top$ is the prefiltering operator, see \cite{marconato2017filter}, performing the first order difference of $h$ and thus the penalty term in (\ref{def_RELS}) penalizes impulses responses for which the norm of the first oder difference of the corresponding $h$ is large
\al{\|g\|_{\Kc_{TC}^{-1}}^2&=(1-\beta)^{-1}\|\Fc^\top h\|^2\nn\\&=(1-\beta)^{-1}\sum_{t=1}^\infty (h_t-h_{t+1})^2.\nn}
The diagonal-correlated (DC) kernel is defined as 
\al{\label{kernelDC} [\Kc_{DC}]_{t,s}  =  \alpha^{|t-s|}\beta^{\max (t,s)}} 
where $0<\beta<1$, $-\beta^{-1/2}<\alpha<\beta^{-1/2}$ and $\eta=[\, \alpha \;\beta \,]^\top$. It is worth noting that we are taking a definition which is not standard, the standard one is $[\Kc_{DC}]_{t,s}  =  \rho^{|t-s|}\beta^{\frac{t+s}{2}}$ and  $\rho=\alpha \beta^{1/2}$, because the former highlights 
the following limits:
\al{\label{condA} \underset{ \alpha \rightarrow 0}{\lim} \Kc_{DC}= \Kc_{DI}, \quad \underset{ \alpha \rightarrow 1}{\lim}\Kc_{DC}=\Kc_{TC} }
that is the DC kernel connects the DI and TC kernel. Indeed, it is not difficult to see that
 \al{\Kc_{DC} = (1-\alpha\beta)(\Fc_\alpha  \Dc \Fc_\alpha^T)^{-1}\nn}
 with 
\al{\label{def_Fa}\Fc_\alpha=\tpl{1,-\alpha,0,\ldots  }.} In plain words, $\alpha$ tunes 
the behavior of the prefiltering operator: $\Fc_\alpha^\top$ behaves as the identity operator for $\alpha$ close to zero, while it behaves as the first order difference operator for $\alpha$ close to one. As a consequence the DC kernel allows to tune the degree of smoothness of $g$.

All these kernels admit a closed form factorization of the inverse and determinant which is an appealing feature for minimizing numerically (\ref{marginal}). Moreover, their inverses are banded matrices: $\Kc_{DI}^{-1}$ is diagonal, $\Kc_{TC}^{-1}$ and $\Kc_{DC}^{-1}$ are tridiagonal.

The stable spline (SS, also called second-order stable spline) kernel induces more smoothness than TC:
\al{ \label{kernelSS}[\Kc_{SS}]_{t,s}  = \frac{\gamma^{t+s}\gamma^{\max (t,s)}}{2}-\frac{\gamma^{3\max (t,s)}}{6}} where $\eta=\gamma$ and $0<\gamma<1$. 
However,  it does not admit a closed form factorization of the inverse and determinant. Moreover, its inverse is not banded. Finally, a kernel with the aforementioned properties which tunes the degree of smoothness and connects TC with SS does not exist.   

 \section{Second-order TC kernel} \label{sec_TC2}
\begin{figure}
\centering
\includegraphics[width=0.5\textwidth]{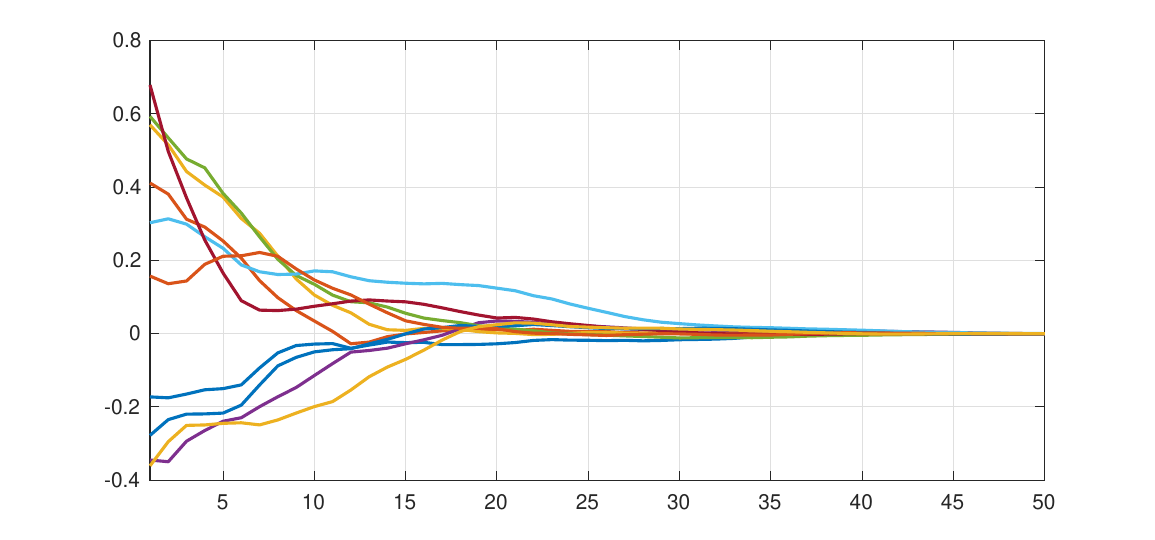}
\caption{Ten realizations of $g\sim \Nc(0, \lambda \Kc_{TC2})$ with $\beta=0.8$ and $\lambda=\|\Kc_{TC2}\|^{-1}$.}
\label{realizationsTC2}
\end{figure}

 In this section we derive a new kernel, hereafter called TC2, which induces more smoothness than TC and represents an alternative to SS. In order to induce more smoothness it is sufficient to take the penalty term as the norm of the second order difference of $h$:
 \al{\|g\|^2_{\Kc^{-1}_{TC2}}=(1-\beta)^{-3}\|(\Fc^T)^2h \|^2,\nn}
 thus 
 \al{\Kc_{TC2}:=(1-\beta)^3(\Fc^2\Dc(\Fc^\top)^2)^{-1}\nn}
 where $\eta=\beta$ and $0<\beta<1$. Figure \ref{realizationsTC2} shows ten realizations of $g$ using the TC2 kernel with $\beta=0.8$. We can notice that the degree of smoothness is similar to the one with $\Kc_{SS}$. 
\smallskip
 \pro \label{pentaTC2}The inverse of $\Kc_{TC2}$ is a pentadiagonal matrix, that is $[(\Kc_{TC2})^{-1}]_{t,s}=0$ for any $|t-s|>3$.\epro
 \smallskip
 \begin{proof} The statement is a particular instance of Proposition \ref{prop_gen_inv}, see Section \ref{sec_high}. \qed
\end{proof}\smallskip\\
Throughout the paper we will use the following result. \smallskip
\lemma[\cite{matriciToeplitz}] \label{lemmaford}Consider a real infinite  lower triangular Toeplitz matrix, defined by the sequence $\{a_k, \; k\geq 0\}$ as follows 
\al{\Xc=\tpl{a_0,a_1,a_2,\ldots}.\nn} If $a_0\neq 0$, $\Xc$ is invertible and the inverse matrix $\Yc=\Xc^{-1}$ is also a lower triangular Toeplitz matrix with elements $\{b_k, \; k\geq 0\}$ given by the following formula
\al{b_0=\frac{1}{a_0}, \; b_k=-\frac{1}{a_0}\sum_{j=0}^{k-1}a_{k-j}b_j\hbox{ for } k\geq 1.\nn}

\elemma
\pro \label{propKTC2} $\Kc_{TC2}$ admits the following closed form expression:  \al{\label{formTC2}[\Kc_{TC2}]_{t,s}=2\beta^{\max(t,s)+1}+(1-\beta)(1+|t-s|)\beta^{\max(t,s)}.}
\epro
\begin{proof} First, $\Fc$ is a lower triangular Toeplitz matrix which is invertible because the main diagonal is composed by strictly positive elements. Therefore, by Lemma \ref{lemmaford} we have 
\al{\Fc^{-2}=\tpl{1,2,\ldots, t,\ldots}.\nn} Moreover,
\al{[\Fc^{-2}]_{t,:}=[\, 0\; \ldots \;0 \;\hspace{-0.5cm} \underbrace{1}_{\tiny \hbox{$t$-th element}} \hspace{-0.5cm}\; 2 \; 3 \ldots \,].\nn}
Therefore,
\al{[\Kc_{TC2}]_{t,s}& =(1-\beta)^3[(\Fc^{-2})^\top \Dc^{-1}\Fc ^{-2}]_{t,s}\nn\\
&=(1-\beta)^3[\Fc^{-2}]_{t,:}\Dc^{-1}[(\Fc^\top)^{-2}]_{:,s}\nn\\
 &=(1-\beta)^3[\Fc^{-2}]_{t,:}\Dc^{-1}[ \Fc^{-2}]_{s,:}^\top\nn\\ 
 &= \sum_{k=\max(t,s)}^\infty \beta^k(k-t+1)(k-s+1).\nn}
Finally, it is not difficult to see that the above series converges to (\ref{formTC2}) by exploiting the identity 
\al{\label{geom_series}\sum_{k=0}^\infty \beta^k=\frac{1}{1-\beta}. }
\qed
\end{proof}\smallskip\\
It is worth noting that the SS kernel is also a second-order generalization of the TC kernel. Indeed, TC and SS are obtained by applying a ``stable'' coordinate change to the first and second order, respectively, spline kernel \cite{PILLONETTO_DENICOLAO2010}. That extension has been 
derived in the continuous time domain, while the one proposed here has been derived in the discrete time domain. 
 
  \section{Second-order DC kernel}\label{sec_DC2}
\begin{figure}
\centering
\includegraphics[width=0.5\textwidth]{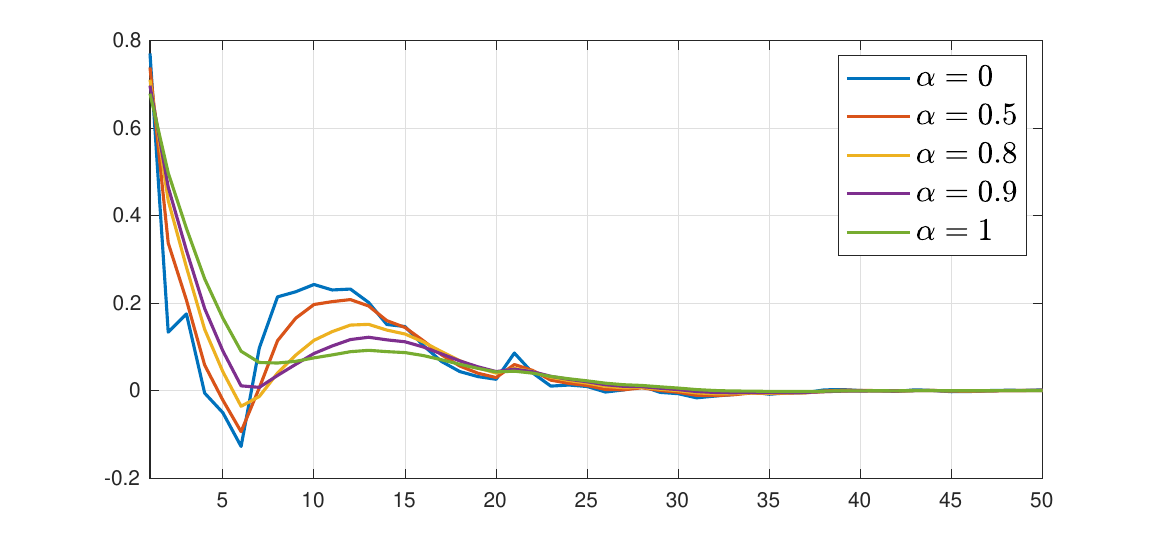}
\caption{One realization of $g\sim \Nc(0,\lambda \Kc_{DC2})$ for different values of $\alpha$. Here, $\lambda=\|\Kc_{DC2}\|^{-1}$.}
\label{fig_transitions}
\end{figure}

  The aim of this section is to introduce a new kernel, hereafter called DC2, which connects the TC and TC2 kernels. The unique difference between TC and TC2 is the prefiltering operator acting on $h$. Thus, the DC2 kernel should perform a transition from $\Fc$ to $\Fc^2$. One possible way is to take
  \al{\label{def_2a}\Fc_{2,\alpha}:=(1-\alpha)\Fc+\alpha \Fc^2 }  
  with $0\leq \alpha\leq 1$ and thus we obtain
  \al{\label{KDC2}\Kc_{DC2}:=\kappa(\Fc_{2,\alpha}\Dc\Fc_{2,\alpha}^\top)^{-1}}
  with $\kappa=(1-\beta)(1-\alpha\beta)(1-\alpha^2\beta)$. In this case we have $\eta=[\, \alpha\; \beta\,]^\top$ with $0<\beta <1$. From the above definition it follows that
  \al{\label{condB}\underset{ \alpha \rightarrow 0}{\lim} \Kc_{DC2}= \Kc_{TC}, \quad \underset{ \alpha \rightarrow 1}{\lim}\Kc_{DC2}=\Kc_{TC2}. } Figure \ref{fig_transitions} shows a realization of the impulse response as a function of $\alpha$ using (\ref{KDC2}); as expected, the degree of smoothness increases as $\alpha$ increases.  \rema It is worth noting that one could consider other transitions, e.g.
  \al{\Fc_{2,\alpha}&=\tpl{1,-1-\alpha^2,\alpha,0,\ldots}\nn\\
  \Fc_{2,\alpha}&=((1-\alpha)\Fc^{-1}+\alpha \Fc^{-2})^{-1}.\nn}
However, as we will see, (\ref{def_2a}) is the unique definition which guarantees that $K_{DC2}$ admits a closed form expression and is the maximum entropy solution of a matrix completion problem.
 \erema
\smallskip
 \pro The inverse of $\Kc_{DC2}$ is a pentadiagonal matrix, that is $[(\Kc_{DC2})^{-1}]_{t,s}=0$ for any $|t-s|>3$.\epro \smallskip
  \begin{proof}
  The proof is similar to the one of Proposition \ref{pentaTC2}.\qed
  \end{proof}\smallskip
  \pro For $0\leq \alpha<1$, $\Kc_{DC2}$ admits the following closed form expression:  {\small \al{\label{formDC2}[\Kc_{DC2}]_{t,s}=\frac{\beta^{\max(t,s)}(1-(1-\beta)\alpha^{|t-s|+1})-\alpha^2\beta^{\max(t,s)+1}}{1-\alpha}.}}
\epro
\begin{proof} First, we notice that
\al{\Fc_{2,\alpha}=((1-\alpha) \Ic + \alpha \Fc) \Fc=\Fc_\alpha \Fc\nn}
where $\Fc_\alpha$ has been defined in (\ref{def_Fa}); $\Ic$ is the identity matrix of infinite dimension. The main diagonal of $\Fc$ and $\Fc_{\alpha}$ is composed by strictly positive elements and thus their inverse exist. By Lemma \ref{lemmaford}, we have 
\al{\Fc^{-1}&=\tpl{1,1, \ldots}\nn\\
\Fc^{-1}_\alpha&=\tpl{1,\alpha,\alpha^2,\ldots}.\nn}
Therefore, 
\al{\Fc_{2,\alpha}^{-1}=\frac{1}{1-\alpha}\tpl{1-\alpha,1-\alpha^2,1-\alpha^3,\ldots}.\nn}
Finally, 
\al{[\Kc_{DC2}&]_{t,s} =\kappa[\Fc_{2,\alpha}^{-1}]_{:,t}^\top \Dc[\Fc_{2,\alpha}^{-1}]_{:,s}\nn\\ &=\kappa\sum_{k=\max(t,s)}^\infty \beta^{k} \frac{1-\alpha^{k-t+1}}{1-\alpha}\frac{1-\alpha^{k-s+1}}{1-\alpha}\nn}
where the above series converges to right hand side of (\ref{formDC2}). The latter fact can be easily proved by using Identity (\ref{geom_series}).
\qed
\end{proof}

\section{Generalized-correlated kernel} \label{sec_GC}
In view of (\ref{condA}) and (\ref{condB}) we can define a general kernel, hereafter called generalized-correlated (GC) kernel, that incorporates the DI, DC, TC, DC2 and TC2 kernels. Let $\Kc_{DI}(\beta)$, $\Kc_{DC}(\alpha,\beta)$, $\Kc_{TC}(\beta)$, $\Kc_{DC2}(\alpha,\beta)$ and $\Kc_{TC2}(\beta)$ be the kernels defined in (\ref{kernelDI}), (\ref{kernelDC}), (\ref{kernelTC}), (\ref{formDC2}) and (\ref{formTC2}), respectively, where we made explicit their dependence on the hyperparameters $0<\alpha<1$   and $0<\beta<1$. Then, we define as GC kernel
\al{\label{kernelGC}\Kc_{GC}(\gamma,\beta)=\left\{\begin{array}{ll}  
  \Kc_{DI}(\beta),&  \gamma=0 \\
   \Kc_{DC}(\gamma,\beta),& 0<\gamma<1  \\
  \Kc_{TC}(\beta),  &  \gamma=1\\
   \Kc_{DC2}(\gamma-1,\beta),  & 1<\gamma<2   \\
  \Kc_{TC2}(\beta),&   \gamma=2.
  \end{array}\right.}
  where $\gamma$ characterizes the smoothness of the impulse response over a wide range. It is worth noting that $\Kc_{GC}$ is a continuous function with respect to $\gamma$ and $\beta$, but not differentiable.

\begin{figure}
\centering
\includegraphics[width=0.5\textwidth]{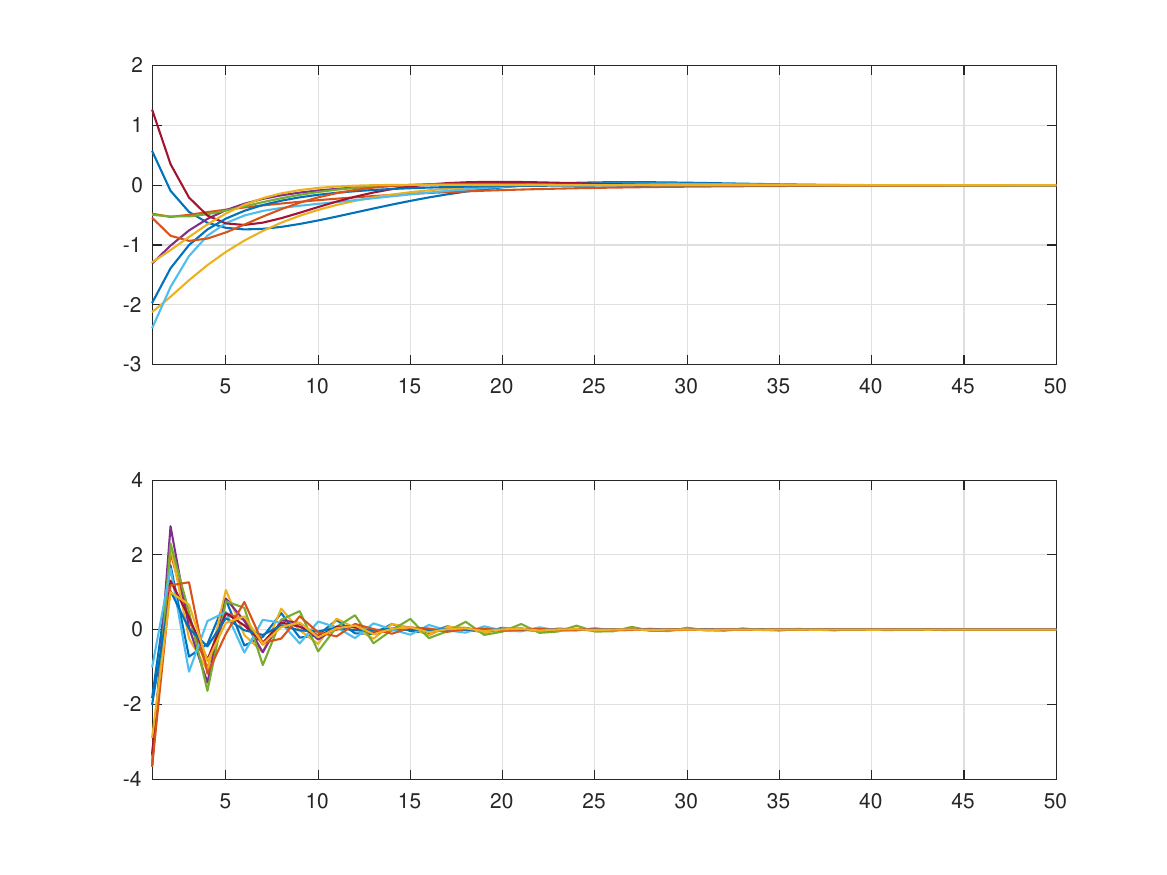}
\caption{{\em Top panel}. Ten realizations of the impulse 
response in the first Monte Carlo study. {\em Bottom panel}. Ten realizations of the impulse 
response in the second Monte Carlo study. }
\label{fig_realsimhl}
\end{figure}

In order to test the superiority of the proposed kernel, in respect to DI, DC, TC, DC2 and TC2, we consider two Monte Carlo studies. The first Monte Carlo study is composed by 200 experiments. In each experiment we generate the impulse response $g$ with practical length $T=50$ as follows:
\al{g_t=\sum_{k=1}^{10} a_k\cos(b_kt+c_k)\nn}
where its parameters are drawn as follows: $a_k\in\Uc([0.2,0.9])$, $b_k\in\Uc([10^{-6}\pi,10^{-1}\pi])$ and $c_k\in\Uc([0,\pi])$. Figure \ref{fig_realsimhl} (top) shows ten realizations drawn from such process. Then, we generate the input of length $N=500$ using the MATLAB function \verb"idinput.m" as a realization drawn from a Gaussian  noise with band [0, 0.6]. Then, we feed the corresponding system (\ref{mod}) with it obtaining the dataset $\mathrm{D}^N:=\{y(t), u(t)\}_{t=1}^N$. Here, $\sigma^2$ is chosen in such a way that the signal to noise ratio is equal to two. Then, we estimate the impulse response using the following estimators:
\begin{itemize}
\item $\hat g_{DI}$ is the estimator in (\ref{def_RELS}) using the diagonal kernel (\ref{kernelDI});
\item $\hat g_{DC}$ is the estimator in (\ref{def_RELS}) using the DC kernel (\ref{kernelDC}); 
\item $\hat g_{TC}$ is the estimator in (\ref{def_RELS}) using the TC kernel (\ref{kernelTC});
\item $\hat g_{D2}$ is the estimator in (\ref{def_RELS}) using the DC2 kernel (\ref{formDC2});
\item $\hat g_{T2}$ is the estimator in (\ref{def_RELS}) using the TC2 kernel (\ref{formTC2});
\item $\hat g_{SS}$ is the estimator in (\ref{def_RELS}) using the SS kernel (\ref{kernelSS});
\item $\hat g_{GC}$ is the estimator in (\ref{def_RELS}) using the GC kernel (\ref{kernelGC}).
\end{itemize}

\begin{figure}[ht]
\centering
\includegraphics[width=0.5\textwidth]{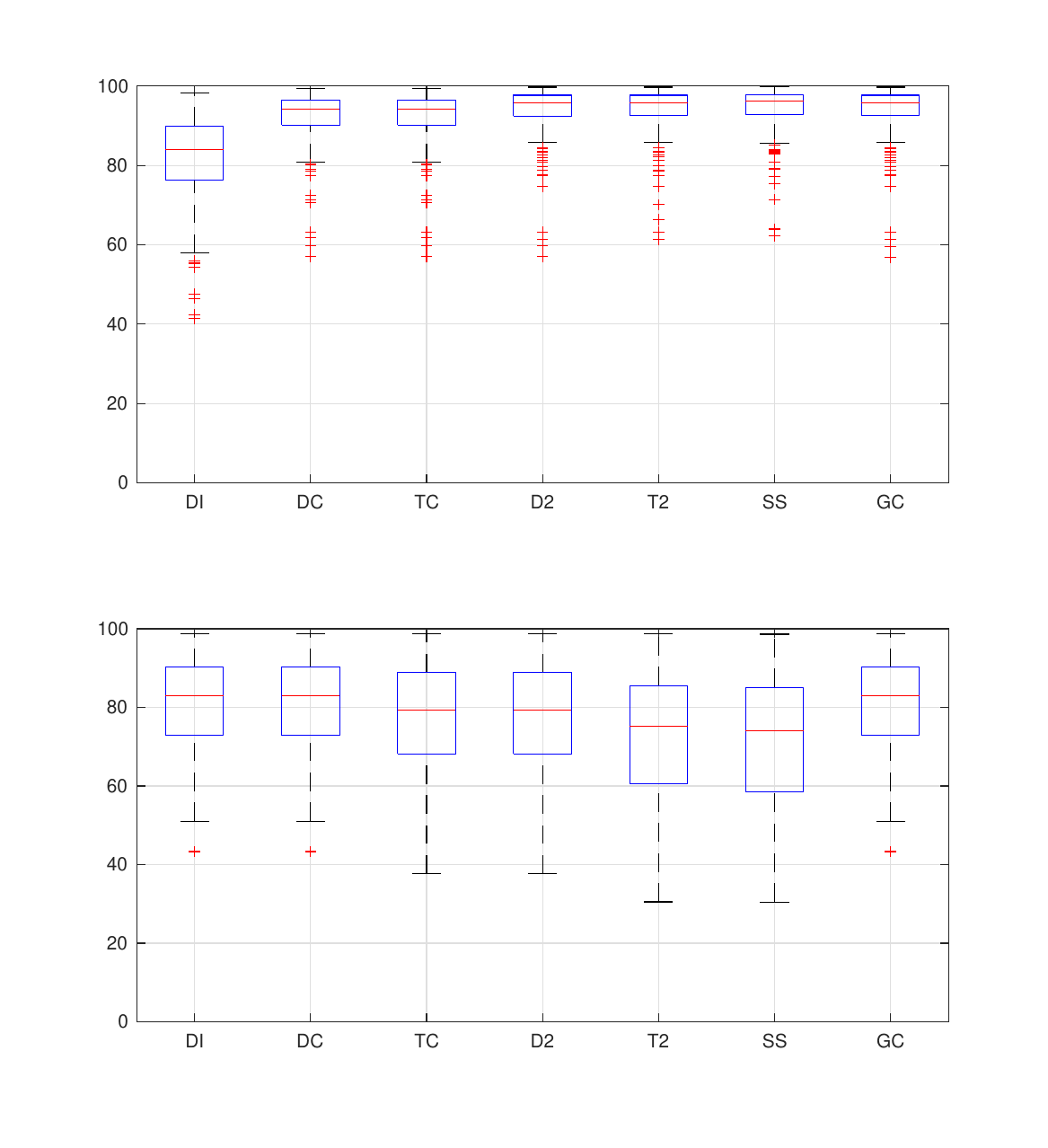}
\caption{Average impulse response fit in the first (top) and second (bottom) Monte Carlo study composed by 200 experiments.}
\label{fig_highlow}
\end{figure}
Finally, for each estimator we compute the average impulse response fit 
\al{\label{def_AIRF}\mathrm{AIRF}=100\left(1- \frac{\|g-\hat g\|}{\|g-\bar g\|} \right)}
where $\bar g=\sum_{t=1}^T g_t$ and $\hat g$ is the corresponding estimator. Clearly, the more $\mathrm{AIRF}$ is close to 100, the better the estimator performance is.   Figure \ref{fig_highlow} (top) shows the boxplot 
of $\mathrm{AIRF}$ for the estimators: D2, T2, SS and GC are the best estimators, while DI is the worst one. In plain words, the best estimators are the ones that are able to induce a sufficient degree of smoothness on the impulse response.

The second Monte Carlo study is likewise to the previous one, but $b_k\in\Uc([0.6\pi,0.7\pi])$. In this case, the realizations of the process $g_t$ are less smooth than before, see Figure \ref{fig_realsimhl} (bottom). Figure \ref{fig_highlow} (bottom) shows the boxplot 
of $\mathrm{AIRF}$ for the estimators: DI, DC and GC are the best estimators, while T2 and SS are the worst ones.  
We conclude that GC is the unique estimator which is able to be well performing in both the situations.

\section{Efficient implementation to estimate the hyperparameters}\label{sec_ME}
 The minimization of (\ref{marginal}) is typically performed through the nonlinear optimization solver \verb"fmincon.m" of Matlab. Thus, the crucial aspect is to consider an efficient algorithm to evaluate (\ref{marginal}). We show that the proposed kernels are suitable for this aim. Recall that $K\in \Sc_T$ denotes the finite dimensional kernel corresponding to  $\Kc$ and defined as 
\al{& [K]_{t,s}=[\Kc]_{t,s},\quad t,s=1\ldots T.\nn}
If $K^{-1}$ admits a closed form expression of its Cholesky factor, then the negative log-marginal likelihood in (\ref{marginal}) can be evaluated efficiently as follows, see \cite{7495008}:  \al{\label{ell_smart}\frac{r^2}{\sigma^2} +(N-T)\log\sigma^2+\log\det(\lambda K)+2\log \det R_1}
where $L$ is the Cholesky factor of $K^{-1}=LL^T$ and $R_1$ is given by the QR factorization 
\al{\left[\begin{array}{cc}R_{d1} & R_{d2} \\ \sigma \sqrt{\lambda^{-1}}L^\top & 0 \end{array}\right]
=QR=Q\left[\begin{array}{cc}
R_1 & R_2  \\0 & r \end{array}\right]\nn}
where $Q^\top Q=I_{T+1}$, $R_1\in\Rs^{T+1\times T}$, $R_2\in\Rs^{T+1}$ and $r\in\Rs$. Moreover, $R_{d1}$ and $R_{d2}$ is given by the QR factorization $[\, A \; y\,]=Q_d[\, R_{d1} \; R_{d2}\,]$ which can be computed ``offline'' before to start the optimization task. In what follows we show that TC2, DC2 and GC admit a closed form expression for $L$ and thus also $\log \det(\lambda K)$. 
\smallskip
\pro \label{prop_TC2_inv_finito}The inverse of $K_{TC2}\in\Sc_T$ admits the following decomposition
\al{K_{TC2}^{-1} =(1-\beta)^{-3}F_T^2 D_T (F_T^2)^\top \nn}
where 
\al{F_T&=\tpl{1,-1,0,\ldots 0}\in\Rs^{T\times T}\nn \\
D_T &=\left[\begin{array}{cc} 	 D_{1,T}  & 0\\ 0& B_T 
 \end{array}\right] \nn\\
 D_{1,T}&=\diag(\beta^{-1},\beta^{-2},\ldots \beta^{T-2})\nn\\
B_T&= (1-\beta)\beta^{-T}\left[\begin{array}{cc}
 \beta+\beta^{2} & 2\beta^{2}	   \\
2\beta^{2} &   1-3\beta +4\beta^{2} \\ 
 \end{array}\right]. \nn} Thus, $K_{TC2}^{-1}$ is a pentadiagonal matrix.
\epro \smallskip
\begin{proof} Consider 
\al{X:=(1-\beta)^3(F_T^2\tilde D_T (F_T^2)^\top)^{-1}\nn} 
where
\al{\label{Dtilda}\tilde D_T=\diag(\beta^{-1},\beta^{-2}\ldots ,\beta^{-T}).}
It is not difficult to see that
\al{F_T^{-2}=\tpl{1,2,\ldots,T}.\nn}
Thus,  by arguments similar to ones used in proof of Proposition \ref{propKTC2}, we have
\al{[X]_{t,s}=(1-\beta)^3\sum_{k=\max(t,s)}^T\beta^k(k-t+1)(k-s+1).\nn}
Without loss of generality, we assume that $t\geq s$; hence,
\al{[X]_{t,s}&=(1-\beta)^3\sum_{k=t}^T\beta^k(k-t+1)(k-s+1).\nn\\
&= (2\beta+(1-b)(1+t-s))\beta^t +\eta(t,s)\nn}
where 
\al{\eta(t,s)=(&1-\beta)(\beta^{T+2}(T-t+1)(T-s+3)\nn\\
&-\beta^{T+1}(T-t+2)(T-s+2))\nn\\ &+2\beta^{T+2}((T-t+1)\beta-(T-t+2))\nn}
where we have exploited the fact that
\al{\sum_{k=0}^T \beta^k =\frac{1-\beta^{T+1}}{1-\beta}.\nn}
Notice that 
\al{[K_{TC2}]_{t,s}=[X]_{t,s}-\eta(t,s)\nn}
and we can rewrite $\eta$ in the shorthand way
\al{\label{structeta}\eta(t,s)=\gamma_1 t+\gamma_1 s+\gamma_2 ts +\gamma_3}
where $\gamma_k$'s are constants not depending on $t$ and $s$. On the other hand, if we take
\al{Y:&=(1-\beta)^3(F_T^2 \Delta^{-1} (F_T^2)^\top)^{-1}\nn\\
&=(1-\beta)^3(F_T^{-2})^\top \Delta F_T^{-2},\nn}
with 
\al{\label{defDelta}\Delta=\left[\begin{array}{ccccc}0  & \ldots & 0 & \ldots & 0 \\\vdots  & \ddots & \vdots &  &\vdots  \\\vdots  &  & 0& z & y \\ 0& \ldots  & 0 & y & x\end{array}\right], } then it is not difficult to see that
\al{[Y]_{t,s}=&(1-\beta)^{3}[-(T(x+2y+z)+x+y)(t+s)\nn\\ &+(x+2y+z)ts+2T(x+y)+x].\nn} By taking into account (\ref{structeta}), we can impose that $x,y,z$ obey the conditions
\al{\gamma_1&=-(1-\beta)^{3}(T(x+2y+z)+x+y)\nn\\
\gamma_2&= (1-\beta)^{3}(x+2y+z)\nn\\
\gamma_3&=(1-\beta)^{3}[2T(x+y)+x].\nn} In this way, $\eta(t,s)=[Y]_{t,s}$. With this choice, we have 
\al{K_{TC2}&=X-Y=(1-\beta)^{3} (F_T^{-2})^\top (\tilde D_T^{-1}-\Delta)(F_T^{-2})\nn\\
&=(1-\beta)^{3} (F_T^{2}(\tilde D_T^{-1}-\Delta)^{-1}(F_T^2)^\top)^{-1}\nn}
where it is not difficult to see that $(\tilde D_T^{-1}-\Delta)^{-1}$ coincides with $D_T$. Finally, the fact that $K_{TC2}^{-1}$ is pentadiagonal follows from Proposition \ref{prop_TC_decomp_Gen}, see Section \ref{sec_high}.\qed
\end{proof}\smallskip
\pro \label{prop_DC2_inv_finito}The inverse of $K_{DC2}\in\Sc_T$ admits the following decomposition
\al{K_{DC2}^{-1} =\kappa^{-1}F_{2,\alpha,T} D_T F_{2,\alpha,T}^\top \nn}
where 
\al{F_{2,\alpha,T}&=(1-\alpha)F_T+\alpha F_T^2\nn \\
D_T &=\left[\begin{array}{cc} 	 D_{1,T}  & 0\\ 0& B_T 
 \end{array}\right] \nn\\
 D_{1,T}&=\diag(\beta^{-1},\beta^{-2},\ldots \beta^{T-2})\nn\\
B_T&= (1-\alpha\beta)\beta^{-T}\nn\\ 
& \times\left[\begin{array}{cc}
 \beta(1+\alpha\beta) & \alpha\beta^{2}(1+\alpha)	   \\
 \alpha\beta^{2}(1+\alpha)&  (1-\beta-\alpha^2\beta) (1-\alpha\beta)+2\alpha^2\beta^2\\ 
 \end{array}\right]. \nn} Thus, $K_{DC2}^{-1}$ is a pentadiagonal matrix.
\epro\smallskip
\begin{proof} Consider 
\al{X:=\kappa(F_{2,\alpha,T}\tilde D_T F_{2,\alpha,T}^\top)^{-1}\nn}
where $\tilde D_T$ has been defined in (\ref{Dtilda}). Notice that $F_{2,\alpha,T}=F_{\alpha,T} F_T$
where 
\al{F_{\alpha,T}=\tpl{1,-\alpha,0 \ldots,0}\in\Rs^{T\times T}\nn}
and 
\al{F_{\alpha,T}^{-1}&=\tpl{1,\alpha^2, \ldots,\alpha^T}\nn\\
F_{2,\alpha,T}^{-1}&=F_T^{-1} F_{\alpha,T}^{-1}\nn\\ &=\frac{1}{1-\alpha}\tpl{1-\alpha,1-\alpha^2, \ldots,1-\alpha^T}\nn.} Without loss of generality, we assume that $t\geq s$, then it is not difficult to see that
\al{[X]_{t,s}&=[(F_T^{-1})^\top(F_{\alpha,T}^{-1})^\top \tilde D_T^{-1}F_{\alpha,T}^{-1}F_T^{-1}]_{t,s}\nn\\
&=\frac{1}{(1-\alpha)^2}\sum_{k=t}^{T}\beta^k (1-\alpha^{k-t+1})(1-\alpha^{k-s+1})\nn\\ 
&= [K_{DC2}]_{t,s}+\eta(t,s)\nn} 
where 
\al{\label{structetaDC} \eta(t,s)=\gamma_1 \alpha^{-t}+\gamma_1\alpha^{-s}+\gamma_2\alpha^{-(t+s)}+\gamma_3}
and $\gamma_k$'s are constants not depending on $t$ and $s$. On the other hand, if we take
\al{Y:&=\kappa(F_{2,\alpha,T} \Delta^{-1} F_{2,\alpha,T}^\top)^{-1}\nn\\
&=\kappa(F_{2,\alpha,T}^{-1})^\top \Delta F_{2\alpha,T}^{-1}\nn}
where $\Delta$ is defined as in (\ref{defDelta}), then it is not difficult to see that 
\al{[Y]_{t,s}=&\kappa [-\alpha^T(z+y+\alpha x+\alpha y )(\alpha^{-t}+\alpha^{-s}) \nn\\
& +\alpha^{2T} (z+2\alpha y+\alpha^2 x) \alpha^{-(t+s)}
+(z+2y+x)].\nn} By taking into account (\ref{structetaDC}), we can impose that $x,y,z$ obey the conditions
\al{\gamma_1&=- \kappa \alpha^T(z+y+\alpha x+\alpha y )\nn\\
\gamma_2&= \kappa \alpha^{2T} (z+2\alpha y+\alpha^2 x)\nn\\
\gamma_3&= \kappa(z+2y+x).\nn
}
In this way, $\eta(t,s)=[Y]_{t,s}$. With this choice, we have 
\al{K_{DC2}&=X-Y=\kappa (F_{2,\alpha,T}^{-1})^\top (\tilde D_T^{-1}-\Delta)F_{2,\alpha,T}^{-1}\nn\\
&=\kappa(F_{2,\alpha,T}(\tilde D_T^{-1}-\Delta)^{-1}F_{2,\alpha,T}^\top)^{-1}\nn}
where it is not difficult to see that $(\tilde D_T^{-1}-\Delta)^{-1}$ coincides with $D_T$. Finally, the fact that $K_{DC2}^{-1}$ is pentadiagonal follows by Proposition \ref{prop_DC_gen_finite} in Section \ref{sec_high}.
\qed
\end{proof} \smallskip\\
By Proposition \ref{prop_TC2_inv_finito} and  \ref{prop_DC2_inv_finito} we have following corollaries.\smallskip
\corr \label{corr1}Let $L$ denote the Cholesky factor of $K_{TC2}^{-1}$, then
\al{[L]_{t,s}=\left\{\begin{array}{lr}\frac{1}{\sqrt{(1-\beta)^3\beta^t }}, & 1 \leq t=s\leq T-2  \\
\frac{-2}{\sqrt{(1-\beta)^3\beta^{t-1} }}, & 2 \leq t=s+1\leq T-1  \\
\frac{1}{\sqrt{(1-\beta)^3\beta^{t-2} }}, & 3 \leq t=s+2\leq T \\
\frac{\sqrt{\beta^{-T+1}(1+\beta)}}{1-\beta}, & t=s=T-1  \\
\frac{-2\sqrt{\beta^{-T+1}}}{(1-\beta)\sqrt{1+\beta}}, & t=s+1=T\\
\sqrt{\frac{\beta^{-T}}{1+\beta}}, & t=s=T \\ 
0,& \hbox{ otherwise.} 
\end{array}\right. \nn} Moreover,
\al{\det K_{TC2}=\beta^{\frac{T(T+1)}{2}}(1-\beta)^{3T-4}.\nn}
\ecorr

\corr Let $L$ denote the Cholesky factor of $K_{DC2}^{-1}$, then
\al{[L]_{t,s}=\left\{\begin{array}{lr}\frac{1}{ \sqrt{\kappa \beta^t }}, & 1 \leq t=s\leq T-2  \\
\frac{-(1+\alpha)}{\sqrt{\kappa \beta^{t-1}}}, & 2 \leq t=s+1\leq T-1 \\
\frac{\alpha}{\sqrt{\kappa \beta^{t-2}}}, & 3 \leq t=s+2\leq T \\
\sqrt{\frac{(1+\alpha \beta)\beta^{-T+1}}{(1-\beta)(1-\alpha^2\beta)}}, & t=s=T-1  \\
\frac{-(1+\alpha)\sqrt{\beta^{-T+1}}}{\sqrt{(1+\alpha \beta ) (1- \beta ) (1-\alpha^2 \beta )}}, & t=s+1=T\\
\sqrt{\frac{\beta^{-T}}{1+\alpha \beta}}, & t=s=T  \\
0,& \hbox{ otherwise.} 
\end{array}\right. \nn} Moreover,
\al{\det K_{DC2}=\beta^{\frac{T(T+1)}{2}}(1-\alpha \beta)^{T-2}(1-\beta)^{T-1}(1-\alpha^2\beta)^{T-1}.\nn}
\ecorr

In view of the above properties, we have 
\al{\log &\det(\lambda K_{TC2})=T\log \lambda+        
\frac{T(T+1)}{2}\log\beta\nn\\ &+(3T-4)\log (1-\beta)\nn\\
\log & \det(\lambda K_{DC2})=T\log \lambda+        
\frac{T(T+1)}{2}\log\beta\nn\\ 
&+(T-2)\log (1-\alpha\beta)+(T-1)\log (1-\beta)\nn\\ &+(T-1)\log (1-\alpha^2\beta).\nn}
In view of the above corollaries and since $K_{DI}^{-1}$, $K_{DC}^{-1}$, $K_{TC}^{-1}$ admit  a closed form expression for the Cholesky factor, see \cite{7495008}, then it follows that the Cholesky factor of $K_{GC}^{-1}$ admits a closed form expression. Accordingly, the minimization of the log-marginal likelihood using GC can be efficiently performed by means of the previous algorithm.
\rema The fact that the inverse kernel matrix is pentadiagonal can be also used to compute efficiently (\ref{def_RELS}) through the alternating direction method of multipliers (ADMM) proposed in \cite{FUJIMOTO21}. Indeed, although that paper considers the case of tridiagonal inverse kernel matrices (e.g. TC and DC kernels)  that idea holds also for banded inverse kernel matrices and the computational flops do not change.\erema

\begin{figure}[ht]
\centering
\includegraphics[width=0.5\textwidth]{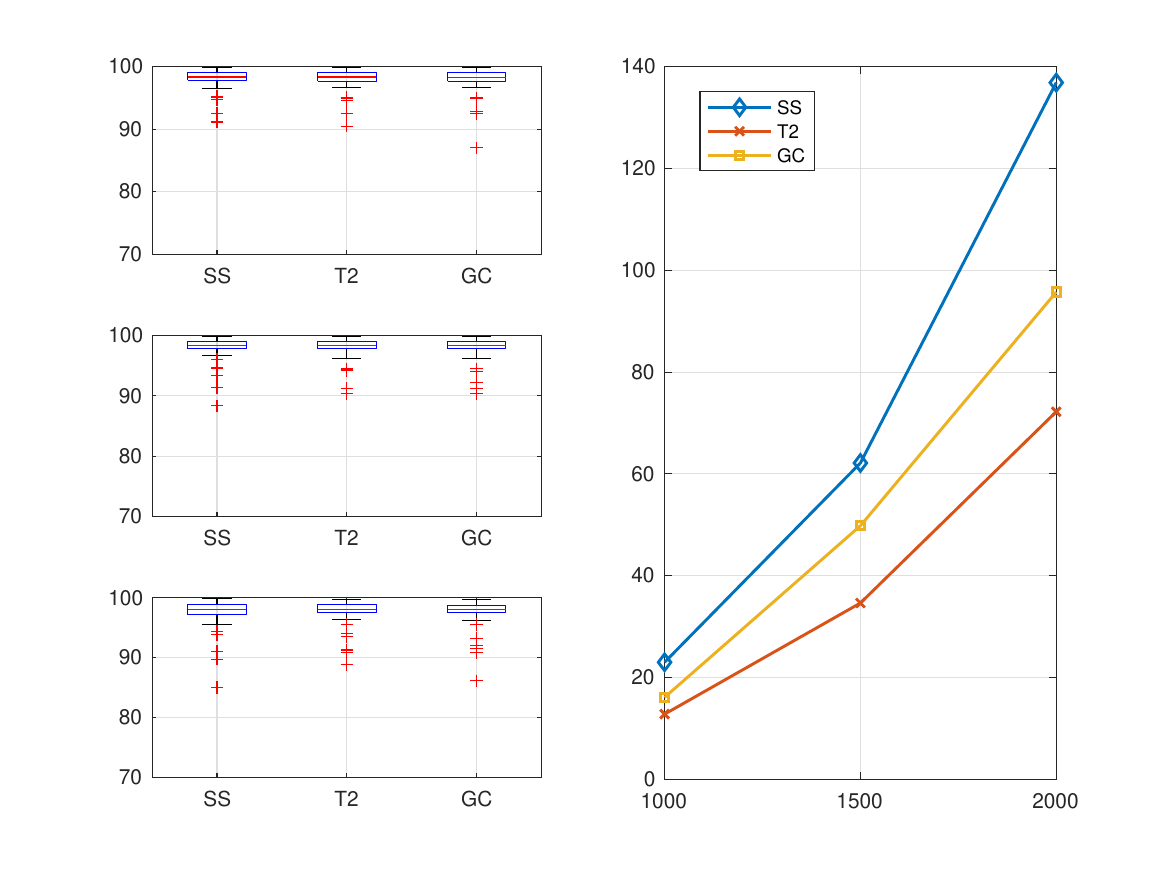}
\caption{{\em Left panels}. Average impulse response fit for $T=1000$ (top), $T=1500$ (middle) and $T=2000$ (bottom). {\em Right panel}.  Average computational time (in seconds) for $T=1000,1500,2000$.}
\label{fig_burden}
\end{figure}

In order to test the aforementioned algorithm equipped with the closed form expressions we consider a Monte Carlo study composed by 50 experiments where the models and  the data are generated likewise to the first Monte Carlo study of Section \ref{sec_GC}, but $a_k\in \mathcal U[0.2, 0.9995]$ and $N=5000$. We consider the following algorithms to estimate the impulse response:
\begin{itemize}
 \item \textbf{T2} is the algorithm in \cite{7495008}, i.e. the one explained before, to compute $\hat g_{T2}$ which exploits the fact that TC2 admits the closed form expression for $L$ and $\log\det(\lambda K)$;
 \item \textbf{GC} is the algorithm in \cite{7495008}, i.e. the one explained before, to compute $\hat g_{GC}$ which exploits the fact that GC admits the closed form expression for $L$ and $\log\det(\lambda K)$;
\item \textbf{SS} is the algorithm in \cite{chen2013implementation} to compute $\hat g_{SS}$ where the Cholesky factor of $K_{SS}$ is computed by \cite[Algorithm 4.2]{doi:10.1137/19M1267349}, i.e. an efficient algorithm  taking linear time which exploits the fact that the SS kernel is extended 2-semiseparable. 
\end{itemize}
For any experiment we measure the computational time (in seconds) of these algorithms through the functions \verb"tic"
and \verb"toc" in Matlab. The simulation is run on a MacBook Air with 3.2GHz Apple M1 processor and 8GB 4266 LPDDR4 memory. Figure \ref{fig_burden} shows the average computational time for the three algorithms using as practical length $T=1000,T=1500,T=2000$ (right panel) and the corresponding average impulse response fit (\ref{def_AIRF}) (left panel).
  While the performance of the estimators is similar, \textbf{T2} exhibits the best computational time and \textbf{SS} the worst one. It is worth noting that the computational time of \textbf{GC} is worse than the one of \textbf{T2} because in the former we have to optimize three hyperparameters (i.e. $\lambda$, $\gamma$ and $\beta$)  while in the latter only two (i.e. $\lambda$ and $\beta$). Finally, for the SS kernel we also considered the algorithm proposed in \cite{7495008} where the Cholesky factor of $K_{SS}$ is computed by \cite[Algorithm 4.2]{doi:10.1137/19M1267349}: the computational time was worse than the one of \textbf{SS}. 
We conclude that SS and TC2 provide a similar performance, thus SS can be safely replaced by TC2 in order to make more efficient the minimization of the negative log-marginal likelihood.

\subsection{Maximum Entropy interpretation}
Proposition \ref{prop_TC2_inv_finito} and Proposition \ref{prop_DC2_inv_finito} are also important to  show that the kernel matrices $K_{TC2}\in\Sc_T$ and $K_{DC2}\in\Sc_T$, with $T\geq 4$, are the maximum entropy solution of a matrix completion problem of the following form. 
\pb[Band extension problem] \label{band_PB}Given $m\in \mathbb N$ and $c_{t,s}$, with $|t-s|\leq m$, find the covariance matrix $\Sigma\in\Sc_T$ of a zero mean Gaussian random vector such that  
\al{ [\Sigma]_{t,s}=c_{t,s}, \quad  |t-s|\leq m.\nn}
\epb
Such an interpretation is important because, as pointed out by Dempster in \cite{dempster1972covariance}, see also \cite{chen2016maximum,chen2018continuous,carli2011maximum,carli2013efficient}, ``the principle of seeking maximum entropy is a principle of seeking maximum simplicity of explanation''. Accordingly, these kernels represent the simplest way of embedding in the prior the fact that the impulse response is BIBO stable and with certain degree of smoothness. Recall that the maximum entropy solution (or extension) of the above problem is defined as 
\al{\label{MEpb} &{\max}_{\Sigma\in\Sc_T } \log \det \Sigma \nn\\
& \hbox{ subject to } [\Sigma]_{t,s}=c_{t,s}, \quad  |t-s|\leq m.}
\smallskip 
\teo \label{teo_1}Consider Problem \ref{band_PB} with $m=2$ and 
\al{c_{t,s}=2\beta^{\max(t,s)+1}+(1-\beta)(1+|t-s|)\beta^{\max(t,s)},\nn} with $|t-s|\leq 2$. Then, the maximun entropy extension solution to (\ref{MEpb}) is $K_{TC2}$.\eteo\smallskip
\begin{proof} In order to prove the statement we need to consider  Problem (\ref{MEpb}) with $m=2\ldots T-2$ where $c_{t,s}=[K_{TC2}]_{t,s}$ with $|t-s|\leq m$. In particular, the solution for $m=2$ is the maximum entropy solution considered in the statement. \smallskip 
\lemma[\cite{dym1981extensions}] 
Problem (\ref{MEpb}) admits solution if and only if 
{\small\al{C_m:=\left[\begin{array}{ccc}c_{t,t}& \ldots & c_{t,m+t} \\ \vdots &  & \vdots \\ c_{t+m,t} & \ldots & c_{t+m,t+m}\end{array}\right]\in\Sc_{m+1},\quad t=1\ldots T-m.\nn}}Under such assumption, the solution is unique with the additional property that its inverse is banded of bandwidth $m$, i.e. its elements in position $(t,s)$ are zero for $|t-s|>m$.  
  \elemma \smallskip 
It is not difficult to see that in our case $C_m\in \Sc_{m+1}$ for $m=2\ldots T-2$ and thus the corresponding band extension problems admit a unique solution. The maximum entropy extension admits a closed form solution that can be computed recursively as follows, see \cite{gohberg1993classes}. Let $\Sigma_T^{(T-2)}$ be the partially specified $T\times   T$ symmetric matrix
\al{\label{defST2} \Sigma_T^{(T-2)}=\left[\begin{array}{ccccc}   
c_{1,1 } & c_{1 ,2 } & \ldots & c_{ 1,T-1 } & x\\
c_{1 ,2 } & c_{ 2, 2} & \ldots & c_{ 2,T-1 } & c_{2 ,T } \\
\vdots & \vdots & & \vdots & \vdots \\ 
c_{1 ,T-1 } & c_{2 ,T-1 } & \ldots & c_{ T-1, T-1} & c_{ T-1, T}\\ 
x & c_{ 2,T } &  \ldots & c_{ T-1,T } & c_{ T, T}\end{array}\right] }
where $x$ is not fixed. Let $X\in \Sc_{T-1}$ be the submatrix of $\Sigma_T^{(T-2)}$ such that 
\al{[X]_{t,s}=c_{t,s}, \; t,s=1\ldots T-1.\nn}
Then, the solution of (\ref{MEpb}) with $m=T-2$, which is called one-step extension, is given by (\ref{defST2})
with 
\al{x= -\frac{1}{y_1} \sum_{j=2}^{T-1} c_{T,j}y_j \nn} 
and $[\, y_1 \; y_2 \; \ldots \; y_{T-1} \,]^\top=L^{-1} [\, 1\; 0\; \ldots \; 0\,]^\top$; moreover, the maximum entropy extension $\Sigma_{ME}$ solution to (\ref{MEpb}) with $2\leq m\leq T-2$ is such that $\Sigma_{ME}^{-1}$ is a band matrix of bandwidth $m$ and for all $m+1<t\leq T$ and $1 \leq s\leq t-m-1$ the submatrix $P^\circ\in\Sc_{t-s+1}$, with $[P^\circ]_{i,j}=[\Sigma_{ME}]_{s-1+i,s-1+j}$, is the one-step extension of the problem
\al{&\min_{P\in\Sc_{t-s+1} } \log \det P \nn\\
& \hbox{ subject to } [P]_{i,j}=c_{s-1+i,s-1+j}, \quad  |i-j|\leq t-s-1.\nn}
Taking into account Proposition \ref{prop_TC2_inv_finito} we know that $K_{TC2}^{-1}$ is banded of bandwidth $m=2$. Let 
\al{\label{defP}P(s,t)=\left[\begin{array}{ccccc}   
c_{s,s } & c_{s ,s+1 } & \ldots & c_{ s,t-1 } & c_{s,t}\\
c_{s ,s+1 } & c_{ s+1, s+1} & \ldots & c_{ s+1,t-1 } & c_{s+1 ,t } \\
\vdots & \vdots & & \vdots & \vdots \\ 
c_{s ,t-1 } & c_{s+1 ,t-1 } & \ldots & c_{ t-1, t-1} & c_{ t-1, t}\\ 
c_{s,t} & c_{ s+1,t } &  \ldots & c_{ t-1,t } & c_{ t, t}\end{array}\right], }   
with $m+1<t\leq T$ and $1 \leq s\leq t-m-1$, be the submatrix of $K_{TC2}$. Then, given the particular definition of $c_{t,s}$'s, it is not difficult to see that
\al{P(s,t-1)&=\beta^{s-1}P(1,t-s)\nn\\ &=\beta^{s-1}(1-\beta)^{3}(F_{t-s}^2 D_{t-s}(F_{t-s}^2)^\top)^{-1} \nn}
where the last equality follows by Proposition \ref{prop_TC2_inv_finito} with $T=t-s$. Then, $P(s,t)$ is the one step-extension of the corresponding band extension problem if $c_{s,t}$ is equal to $x$ and the latter is given as follows. We define
\al{[\, & y_1 \; y_2 \; \ldots \; y_{t-s} \,]^\top=P(s,t-1)^{-1}[\, 1\; 0\; \ldots \; 0\,]^\top\nn\\ &=\beta^{1-s}(1-\beta)^{-3}F_{t-s}^2 D_{t-s}(F_{t-s}^2)^\top [\, 1\; 0\; \ldots \; 0\,]^\top\nn\\ 
& =\beta^{1-s}(1-\beta)^{-3} [\, \beta^{-1}\; -2\beta^{-1}\; \beta^{-1}\; 0\ldots \; 0\,]^\top;\nn} therefore 
\al{x&=-y_1^{-1} \sum_{j=2}^{t-s} c_{t,s+j}y_j\nn\\
&=-\beta(-2\beta^{-1}c_{s+1,t}+\beta^{-1}c_{s+2,t})=2c_{s+1,t}-c_{s+2,t}\nn\\
&=2\beta^{t+1}+(1-	\beta)(1+t-s)\beta^t=[K_{TC2}]_{t,s}\nn}
which concludes the proof.\qed
\end{proof}\smallskip\\
\teo Consider Problem \ref{band_PB} with $m=2$ and 
\al{c_{t,s}=\frac{\beta^{\max(t,s)}(1-(1-\beta)\alpha^{|t-s|+1})-\alpha^2\beta^{\max(t,s)+1}}{1-\alpha},\nn} with $|t-s|\leq 2$. Then, the maximun entropy extension solution to (\ref{MEpb}) is $K_{DC2}$.\eteo\smallskip
\begin{proof} The proof is similar to the one of Theorem \ref{teo_1}. More precisely, in this case we have $P(t,s)$ is the submatrix of $K_{DC2}$  defined as in (\ref{defP})
with $m+1<t\leq T$ and $1 \leq s\leq t-m-1$. Then, given the particular definition of $c_{t,s}$'s, it is not difficult to see that
\al{P(s,t-1)&=\beta^{s-1}P(1,t-s)\nn\\ &=\beta^{s-1}\kappa(F_{2,\alpha,t-s}D_{t-s}F_{2,\alpha,t-s}^\top)^{-1} \nn}
where the last equality follows by Proposition \ref{prop_DC2_inv_finito} with $T=t-s$. Then, $P(s,t)$ is the one step-extension of the corresponding band extension problem if $c_{s,t}$ is equal to $x$ and the latter is given as follows. We define
\al{[\, & y_1 \; y_2 \; \ldots \; y_{t-s} \,]^\top\nn\\ &=\beta^{1-s}\kappa^{-1}F_{2,\alpha,t-s} D_{t-s}F_{2,\alpha,t-s}^\top [\, 1\; 0\; \ldots \; 0\,]^\top\nn\\ 
& =\beta^{1-s}\kappa^{-1} [\, \beta^{-1}\; -(1+\alpha)\beta^{-1}\; \alpha\beta^{-1}\; 0\ldots \; 0\,]^\top;\nn} therefore 
\al{x&=-y_1^{-1} \sum_{j=2}^{t-s} c_{t,s+j}y_j \nn\\
&=-\beta(-(1+\alpha)\beta^{-1}c_{s+1,T}+\alpha\beta^{-1}c_{s+2,T})\nn\\  &=(1+\alpha)c_{s+1,t}-\alpha c_{s+2,t}\nn\\
&=\frac{\beta^{t}(1-(1-\beta)\alpha^{t-s+1})-\alpha^2\beta^{t}}{1-\alpha}=[K_{DC2}]_{t,s}\nn }
which concludes the proof.\qed
\end{proof}

 \section{Higher-order extensions}\label{sec_high}
\begin{figure}
\centering
\includegraphics[width=0.5\textwidth]{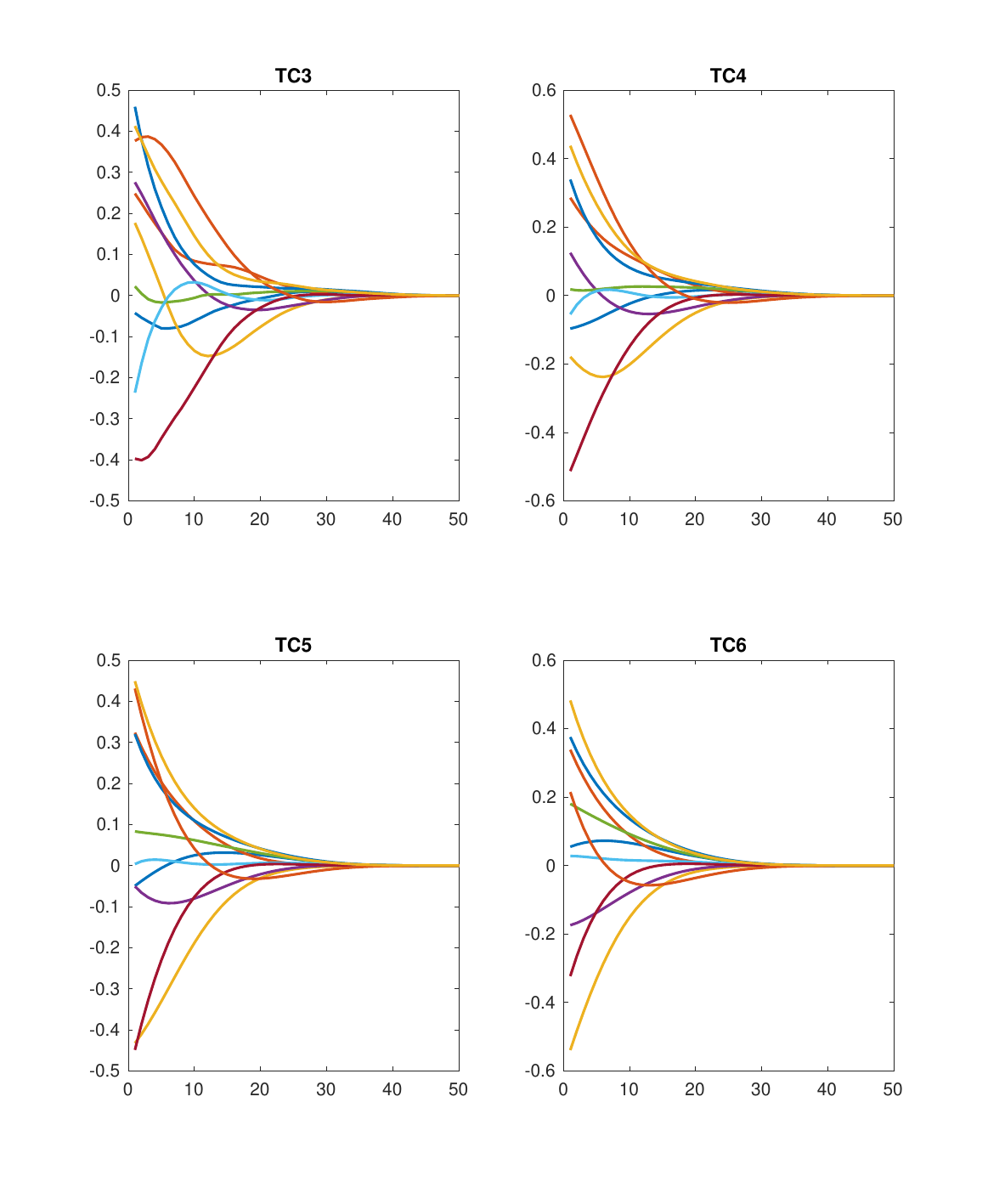}
\caption{Ten realizations of $g\sim \Nc(0, \lambda \Kc_{TC\delta})$ with $\beta=0.8$, $\delta=3,4,5,6$ and $\lambda=\|\Kc_{TC\delta}\|^{-1}$.}
\label{realizationsTCdelta}
\end{figure}

Drawing inspiration from Section  \ref{sec_TC2} 
we can define the TC kernel of order $\delta\in \Ns$ as  
\al{\label{kernelTCd} \Kc_{TC\delta}=\kappa_\delta(\Fc^\delta \Dc (\Fc^\delta)^\top)^{-1}}
where $\kappa_\delta$ is a suitable normalization constant. Here, $\eta=\beta$ with $0<\beta<1$. Figure \ref{realizationsTCdelta} shows ten realizations of $g$ using the TC$\delta$ kernel with $\beta=0.8$ and for different values of $\delta$.    As expected, the larger $\delta$ is the more smoothness is induced on $g$. \smallskip

\pro \label{prop_gen_inv}The inverse of $\Kc_{TC\delta}$ is a banded matrix of bandwidth $\delta$, that is $[\Kc_{TC\delta}^{-1}]_{t,s}=0$ for any $|t-s|>\delta$.\epro\smallskip
 \begin{proof} We prove the claim by induction. First, for $\delta=1$ we have that TC$\delta$ is the standard TC and its inverse is tridiagonal, i.e. the claim holds. Assume that $\Kc_{TC\delta-1}^{-1}$ is a banded matrix of bandwidth $\delta-1$. Then,
 \al{\Kc_{TC\delta}^{-1}= \frac{\kappa_{\delta-1}}{\kappa_{\delta}} \Fc\Kc_{TC\delta-1}^{-1}\Fc^\top.\nn}
Notice that $\Fc=\Ic-\Sc$ where $\Sc$ is the lower shift matrix and $\Ic$ the identity matrix, both infinite dimensional. Hence,
\al{\label{decomp_shift}&\Kc_{TC\delta}^{-1}=\kappa_{\delta-1}  \kappa_{\delta}^{-1} [\Kc_{TC\delta-1}^{-1}\nn\\ &\hspace{0.2cm}+\Sc\Kc_{TC\delta-1}^{-1}\Sc^\top-\Kc_{TC\delta-1}^{-1}\Sc^\top-\Sc\Kc_{TC\delta-1}^{-1}].}  It is well known that premultiplying a matrix $A$ by a lower shift matrix results in the elements of $A$ being shifted downward by one position, with zeroes appearing in the top row. Thus, in view of (\ref{decomp_shift}), we have that 
$\Sc\Kc_{TC\delta-1}^{-1}\Sc^\top$ is a band matrix with bandwidth $\delta-1$, while $\Kc_{TC\delta-1}^{-1}\Sc^\top+\Sc\Kc_{TC\delta-1}^{-1}$ and thus $\Kc_{TC\delta-1}$ are band matrices with bandwidth $\delta$. \qed
\end{proof} \smallskip

Also in this case one could try to find the closed form expression for $\Kc_{TC\delta}$, however its derivation is not straightforward from the case $\delta=2$. On the other hand, we can define the corresponding finite dimensional kernel matrix $K_{TC\delta}\in \Sc_T$ as
\al{[K_{TC\delta}]_{t,s}=[\Kc_{TC\delta}]_{t,s}, \; \; t,s=1\ldots T.\nn}

\pro \label{prop_TC_decomp_Gen}The finite dimensional kernel $K_{TC\delta}$ admits the following decomposition: \al{K_{TC\delta}^{-1}=\kappa_\delta^{-1}F_T^\delta D_T (F_T^\delta)^\top\nn}
 where 
 \al{D_T &=\left[\begin{array}{cc} 	 D_{1,T}  & 0\\ 0& B_T 
 \end{array}\right] \nn\\
 D_{1,T}&=\diag(\beta^{-1},\beta^{-2},\ldots \beta^{T-\delta}),
 \nn} and $B_T$ is a $\delta\times \delta$ matrix. Thus,  $K_{TC\delta}^{-1}$ is banded of bandwidth $\delta$.
\epro\smallskip
\begin{proof}
Let $\Vc^{(j)}\in\Rs^{\infty \times T}$ denote a matrix whose first $j-1$ columns coincide with the null sequence and the remaining ones do not, thus $\Vc^{(T+1)}$ is the null matrix. We use $\sim$ to denote the equivalence relation $\Xc\sim \Yc$ which means that $\Xc\in\Rs^{\infty \times T}$ and $\Yc\in\Rs^{\infty \times T}$ have the first columns (in the same number) equal to the null sequence and the other ones do not. Thus, the latter induces a splitting of $\Rs^{\infty\times T}$ through the corresponding equivalence classes $[ \Vc^{(j)}]=\{ \Xc\in \Rs^{\infty\times T} \hbox{ s.t. } \Xc\sim \Vc^{(j)}\}$ with $1\leq j\leq T+1$. In what follows, in order to ease the exposition (and thus with some abuse of notation) we use the symbol $=$ instead of $\sim$ in all the (submatrix) relations involving $\Vc^{(j)}$, with $j=1\ldots T+1$.

First, notice that $\Fc=\Ic-\Sc$ and $F_T=I_T-S$ where $\Sc$
 and $S$ denote, respectively, the infinite and finite dimensional lower shift matrix.  Recall that postmultiplying $\Vc^{(j)}$, with $1\leq j \leq T$, by $S$ results in the columns of  $\Vc^{(j)}$ being shifted left by one position with a null sequence appearing in the last column position, thus
 \al{\label{eq1shift} \Vc^{(j)}F_T=\Vc^{(j-1)};}
premultiplying $\Vc^{(j-1)}$, with $1\leq j \leq T$, by $\Sc$
results in the rows of  $\Vc^{(j-1)}$ being shifted downward  by one position with a null row vector appearing in the first top row, thus
 \al{\label{eq2shift} \Vc^{(j-1)}F_T=\Vc^{(j-1)}.}
Combining (\ref{eq1shift})-(\ref{eq2shift}), we obtain 
\al{ \Vc^{(j)}F_T=\Fc\Vc^{(j-1)}\nn}
and thus
\al{\label{eq3shift} \Fc^{-1}\Vc^{(j)}F_T=\Vc^{(j-1)}, \; \; 1\leq j\leq T.}
Then, we have
\al{\label{eq4shift}\Fc^{-1} \left[\begin{array}{c} 	 I_T \\ \Vc^{(j)} 
 \end{array}\right]F_T=\left[\begin{array}{c} 	 I_T \\ \Oc+ 
 \Fc^{-1}\Vc^{(j)}F_T\end{array}\right] }
where $\Oc\in\Rs^{\infty \times T}$ is a matrix whose last column is a sequence of ones, while the other columns are null sequencess, i.e. $\Oc=\Vc^{(T-1)}$. Accordingly, by (\ref{eq3shift})-(\ref{eq4shift}) we have 
\al{\label{eq5shift}\Fc^{-1} \left[\begin{array}{c} 	 I_T \\ \Vc^{(j)} 
 \end{array}\right]F_T=\left[\begin{array}{c} 	 I_T \\  \Vc^{(j-1)} \end{array}\right], \; \; 1\leq j\leq T+1. }
Notice that
\al{K_{TC\delta}&= \left[\begin{array}{cc} 	 I_T  & 0
 \end{array}\right] \Kc_{TC\delta}\left[\begin{array}{c} 	 I_T \\0 
 \end{array}\right]  \nn\\ 
 &=\kappa_\delta\left[\begin{array}{cc} 	 I_T  & 0
 \end{array}\right](\Fc^{-\delta})^\top\Dc^{-1}\Fc^{-\delta } \left[\begin{array}{c} 	 I_T \\0 
 \end{array}\right].\nn
 } 
Consider
\al{Y:&= (F_T^{\delta})^\top K_{TC\delta}F_T^{\delta} \nn\\
&=\kappa_\delta \Wc_\delta^\top\Dc^{-1} \Wc_\delta\nn}
where 
\al{\Wc_\delta=\Fc^{-\delta } \left[\begin{array}{c} 	 F_T^\delta \\0 
 \end{array}\right]=\Fc^{-\delta } \left[\begin{array}{c} 	 F_T^\delta \\ \Vc^{(T+1)}
 \end{array}\right].\nn}
Then, it remains to prove that $Y=\kappa_\delta D_T^{-1}$. Indeed, 
\al{\Wc_\delta&=\Fc^{-(\delta-1) } \Fc^{-1}\left[\begin{array}{c} 	 I_T\\ \Vc^{(T+1)} 
 \end{array}\right]F_T F_T^{\delta-1}\nn\\ 
 &=\Fc^{-(\delta-1) }\left[\begin{array}{ccc} I_T \\ \Vc^{(T)}\end{array}\right]F_T^{\delta-1}\nn\\ 
 &= \ldots =\left[\begin{array}{ccc} I_T \\ \Vc^{(T+1-\delta)}\end{array}\right]\nn}
 where we exploited (\ref{eq5shift}). Thus, 
 \al{Y&=\kappa_\delta\left[\begin{array}{cc} I_T & (\Vc^{(T+1-\delta)})^\top\end{array}\right] \Dc^{-1}\left[\begin{array}{ccc} I_T \\ \Vc^{(T+1-\delta)}\end{array}\right]\nn\\
 &=\kappa_\delta\left[\begin{array}{cc} I_T & (\Vc^{(T+1-\delta)})^\top\end{array}\right]
 \left[\begin{array}{cc} D_{1,T}^{-1} & 0\\ 
 0& \tilde \Dc^{-1} \end{array}\right] 
 \left[\begin{array}{ccc} I_T \\ \Vc^{(T+1-\delta)}\end{array}\right]\nn\\
 &=\kappa_\delta (D_{1,T}^{-1}+(\Vc^{(T+1-\delta)})^\top \tilde \Dc^{-1} \Vc^{(T+1-\delta)})=\kappa_\delta D_T^{-1}\nn} 
where $\tilde \Dc=\diag(\beta^{T-\delta-1},\beta^{T-\delta-2},\ldots)$.
\qed
\end{proof}\smallskip\\


It remains to design the DC kernel of oder $\delta$ connecting $\Kc_{TC\delta-1}$ and $\Kc_{TC	\delta}$. Drawing inspiration from Section \ref{sec_DC2} we define it as 
\al{\label{KDCd}\Kc_{DC\delta}=\kappa_{\delta} (\Fc_{\delta,\alpha} \Dc \Fc_{\delta,\alpha}^\top)^{-1}}
where 
  \al{\Fc_{\delta,\alpha}:=(1-\alpha)\Fc^{\delta-1}+\alpha \Fc^\delta\nn }
  and $\kappa_\delta$ is the normalization constant. Here, $\eta=[\, \beta\; \alpha\,]^\top$ with $0<\beta<1$ and $0\leq \alpha\leq1$.
In Figure \ref{fig_transitions} we show a realization of the impulse response using (\ref{KDCd}) with $\delta=3$ as a function of $\alpha$; as expected, the degree of smoothness increases as $\alpha$ increases.\smallskip
\begin{figure}
\centering
\includegraphics[width=0.5\textwidth]{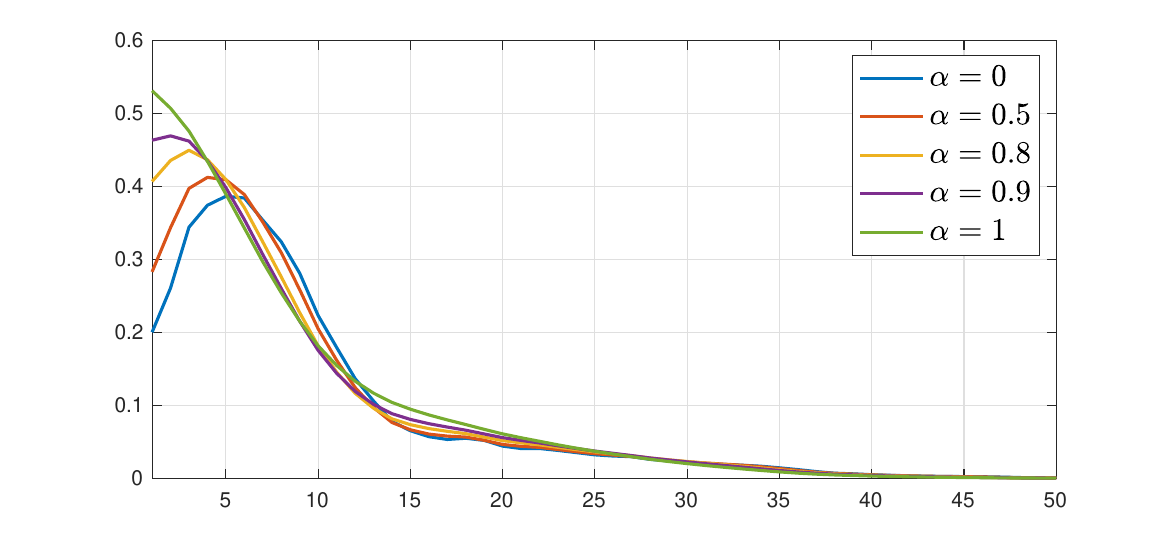}
\caption{One realization of $g\sim \Nc(0,\lambda \Kc_{DC3})$ for different values of $\alpha$. Here, $\beta=0.8$ and $\lambda=\|\Kc_{DC3}\|^{-1}$.}
\label{trans23}
\end{figure}
\pro \label{prop_gen_invDC}The inverse of $\Kc_{DC\delta}$ is a banded matrix of bandwidth $\delta$, that is $[\Kc_{DC\delta}^{-1}]_{t,s}=0$ for any $|t-s|>\delta$.\epro\smallskip
\begin{proof}
First, for $\delta=1$ $K_{DC\delta}$ is the standard DC kernel whose inverse is tridiagonal, i.e. the statement holds.
 Finally, notice that
\al{\Fc_{\delta,\alpha}:= \Fc((1-\alpha)\Fc^{\delta-2}+\alpha \Fc^{\delta-1})=\Fc\Fc_{\delta-1,\alpha},\nn }
thus
\al{\Kc_{DC\delta}^{-1}=\kappa_{\delta-1}\kappa_{\delta}^{-1}\Fc \Kc_{DC\delta-1}^{-1}\Fc^\top.\nn}
 Accordingly, the remaining part of the proof is similar to the one of Proposition \ref{prop_gen_inv}.\qed 
\end{proof}\smallskip\\

Also in this case the finite dimensional kernel $K_{DC\delta}\in\Sc_T$ is defined as 
\al{[K_{DC\delta}]_{t,s}= [\Kc_{DC\delta}]_{t,s}, \; \; t,s=1\ldots T.\nn}
\pro\label{prop_DC_gen_finite} The finite dimensional kernel $K_{DC\delta}$ admits the following decomposition:  
\al{K_{DC\delta}^{-1}=\kappa_{\delta} F_{\delta,\alpha,T} D_T (F_{\delta,\alpha,T})^\top\nn}
 where 
\al{F_{\delta,\alpha,T}&=(1-\alpha)F_T^{\delta-1}+\alpha F_T^\delta\nn \\
D_T &=\left[\begin{array}{cc} 	 D_{1,T}  & 0\\ 0& B_T 
 \end{array}\right] \nn\\
 D_{1,T}&=\diag(\beta^{-1},\beta^{-2},\ldots \beta^{T-\delta}) \nn}   and $B_T$ is a $\delta\times \delta$ matrix; Thus,  $K_{DC\delta}^{-1}$ is banded of bandwidth $\delta$.
\epro\smallskip
\begin{proof}
The proof is similar to the one of Proposition \ref{prop_TC_decomp_Gen}.\qed
\end{proof}\smallskip
\\
Finally, this extension can be applied also to the high-frequency  (HF) kernel, see \cite{6160606}: \al{[\Kc_{HF}]_{t,s} =  (-1)^{ |t-s|} \beta^{\max(t,s)}=(-1)^{ |t-s|}[\Kc_{TC}]_{t,s}  
\nn } where $0<\beta<1$. We define the high frequency kernel of oder $\delta\in\Ns$ as
\al{[\Kc_{HF\delta}]_{t,s}=(-1)^{|t-s|} [\Kc_{TC\delta}]_{t,s}.\nn}
Moreover, we can define the high frequency diagonal-correlated (HC) kernel connecting HF$\delta-1$ and HF$\delta$ as 
\al{[\Kc_{HC}]_{t,s}=(-1)^{|t-s|} [\Kc_{DC\delta}]_{t,s}.\nn}
It is straightforward to see that $\Kc_{HF\delta}^{-1}$ and $\Kc_{HC\delta}^{-1}$ are banded of bandwidth $\delta$, as well as their finite dimensional matrices $K_{HF\delta}^{-1}$ and $K_{HC\delta}^{-1}$. It is possible to find the closed form expression for the Cholesky factor and the determinant of  $K_{HF2}^{-1}$ and $K_{HC2}^{-1}$. Finally, $K_{HF2}$ and $K_{HC2}$ are, respectively, the maximum entropy solution of a band extension problem similar to the ones introduced in Section \ref{sec_ME}.

\section{Frequency analysis}\label{sec_freq}
An exponentially convex local stationary (ECLS) kernel $\Kc\in\Sc_\infty$ admits the following decomposition
\al{\label{defECLS}[\Kc]_{t,s}=\beta^{\frac{t+s}{2}} [\Wc]_{t,s}}
where $\Wc\in\Sc_\infty$ is a stationary kernel, i.e. the covariance function of a stationary process and thus $[\Wc]_{t,s}=[\Wc]_{t+k,s+k}$ for any $k\in\Ns$. Recall that TC, DC and SS are ECLS kernels.
It is straightforward to see that TC2 and DC2 are ECLS kernel 
whose stationary parts are, respectively,
\al{[\Wc_{TC2}]_{t,s}&= 2\beta^{\frac{|t-s|}{2}+1}+(1-\beta)(1+|t-s|)\beta^{\frac{|t-s|}{2}}\nn\\
[\Wc_{DC2}]_{t,s}&=\frac{\beta^{\frac{|t-s|}{2}} (1-(1-\beta)\alpha^{|t-s|+1})-\alpha^2 \beta^{\frac{|t-s|}{2}+1}}{1-\alpha}. \nn}   
\teo TC$\delta$ and DC$\delta$ kernels with $\delta>2$ are ECLS, that is
\al{ \Kc_{TC\delta}=\beta^{\frac{t+s}{2}} [\Wc_{TC\delta}]_{t,s}, \; \;\Kc_{DC\delta}=\beta^{\frac{t+s}{2}} [\Wc_{DC\delta}]_{t,s} \nn}
where $\Wc_{TC\delta}$ and $\Wc_{DC\delta}$ are stationary kernels.
\eteo \smallskip\begin{proof} We only prove the claim for TC$\delta$ because the one for DC$\delta$ is similar. By (\ref{kernelTCd}), we have that 
\al{\label{def_TCd_ecls}\Kc_{TC\delta}=\kappa_\delta \Xc^\top \Dc^{-1}\Xc}
where $\Xc=\Fc^{-\delta}=(\Fc^{-1})^\delta$. Since $\Fc$ is lower triangular, Toeplitz and invertible, then by Lemma \ref{lemmaford} we know that $\Fc^{-1}$ is lower triangular and Toeplitz. Accordingly, $\Xc$ is lower triangular and Toeplitz because it is given by a product of lower triangular and Toeplitz matrices. Hence, let 
\al{\Xc=\tpl{x_1,x_2,x_3,\ldots}. \nn} 
Moreover,
\al{[\Xc]_{t,:}=[\, 0\; \ldots \;0 \;\hspace{-0.5cm} \underbrace{x_1}_{\tiny \hbox{$t$-th element}} \hspace{-0.5cm}\; x_2 \; x_3 \ldots \,].\nn}
Taking into account (\ref{def_TCd_ecls}), we have 
\al{\label{TCdECLS}[\Kc_{TC\delta}]_{t,s}&=\kappa_\delta [\Xc]_{t,:}\Dc^{-1} [\Xc]_{s,:}^\top\nn\\ 
&=\kappa_\delta   \sum_{k=1}^\infty \beta^{\max(t,s)+k-1} x_k x_{k+|t-s|}\nn\\
&=\beta^{\frac{t+s}{2}}\underbrace{ \kappa_\delta  \sum_{k=1}^\infty \beta^{\frac{|t-s|}{2}+k-1} x_k x_{k+|t-s|}}_{:=[\Wc_{TC\delta}]_{t,s}}}
where we have exploited the fact that $\max(t,s)=(t+s)/2+|t-s|/2$. It is straightforward to see that $\Wc_{TC\delta}$ is a stationary kernel. 
In view of (\ref{defECLS}) and (\ref{TCdECLS}), we conclude that TC$\delta$ is ECLS.
\qed \end{proof}\smallskip

Although it is not immediate to derive the closed form expression for $W_{TC\delta}$ and $W_{DC\delta}$, we can compute them numerically:
\al{[W_{TC\delta}]_{t,s}\approx \beta^{-\frac{t+s}{2}\nn}[K_{TC\delta}]_{t,s}}
and likewise for DC$\delta$. Clearly, the larger $T$ is, the better the approximation above is.  

\begin{figure}
\centering
\includegraphics[width=0.5\textwidth]{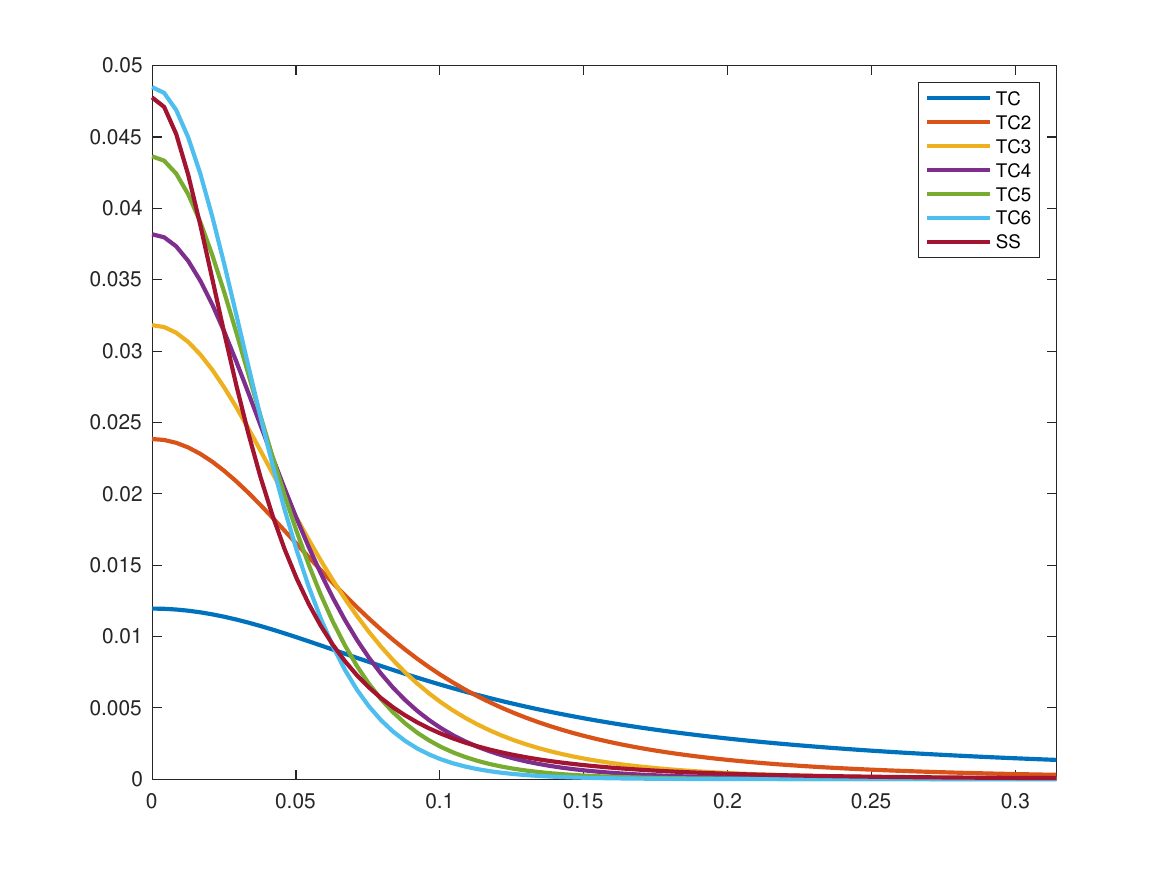}
\caption{Power spectral density of the stationary part of TC, TC$\delta$, with $\delta=2\ldots 6$, and SS with $\beta=0.8$. All those power spectral densities are normalized to one in order to ease the comparison.}
\label{fig_freq}
\end{figure}

Therefore, it is interesting to compare the frequency content in their stationary parts. In doing that, we 
recall that 
\al{[W]_{t,s}=\frac{1}{2\pi}\int_{-\pi}^\pi \phi(\vartheta)\cos(\vartheta(t-s))\d \vartheta\nn} 
where $\phi(\vartheta)$, with $\vartheta\in[0,2\pi]$, is the power spectral density of the (stationary) process. In order to compare SS with the others we need to choose $\gamma=\sqrt[3]{\beta}$ in (\ref{kernelSS}); in this way the latter has the exponential part as in (\ref{defECLS}). Figure \ref{fig_freq} shows the 
  power spectral densities of the stationary part of TC, TC$\delta$, with $\delta=2\ldots 6$ and SS. As expected, the higher TC$\delta$ is, the more statistical power is concentrated for  frequencies close to zero. TC2 promote less smoothness than SS, while the latter is more similar to TC5 and TC6. It is worth noting that we can plot also the power spectral density corresponding to DC$\delta$. The latter smoothly changes from the one of TC$\delta-1$, with $\alpha=0$, to the one of TC$\delta$, with $\alpha=1$.

\begin{figure}
\centering
\includegraphics[width=0.5\textwidth]{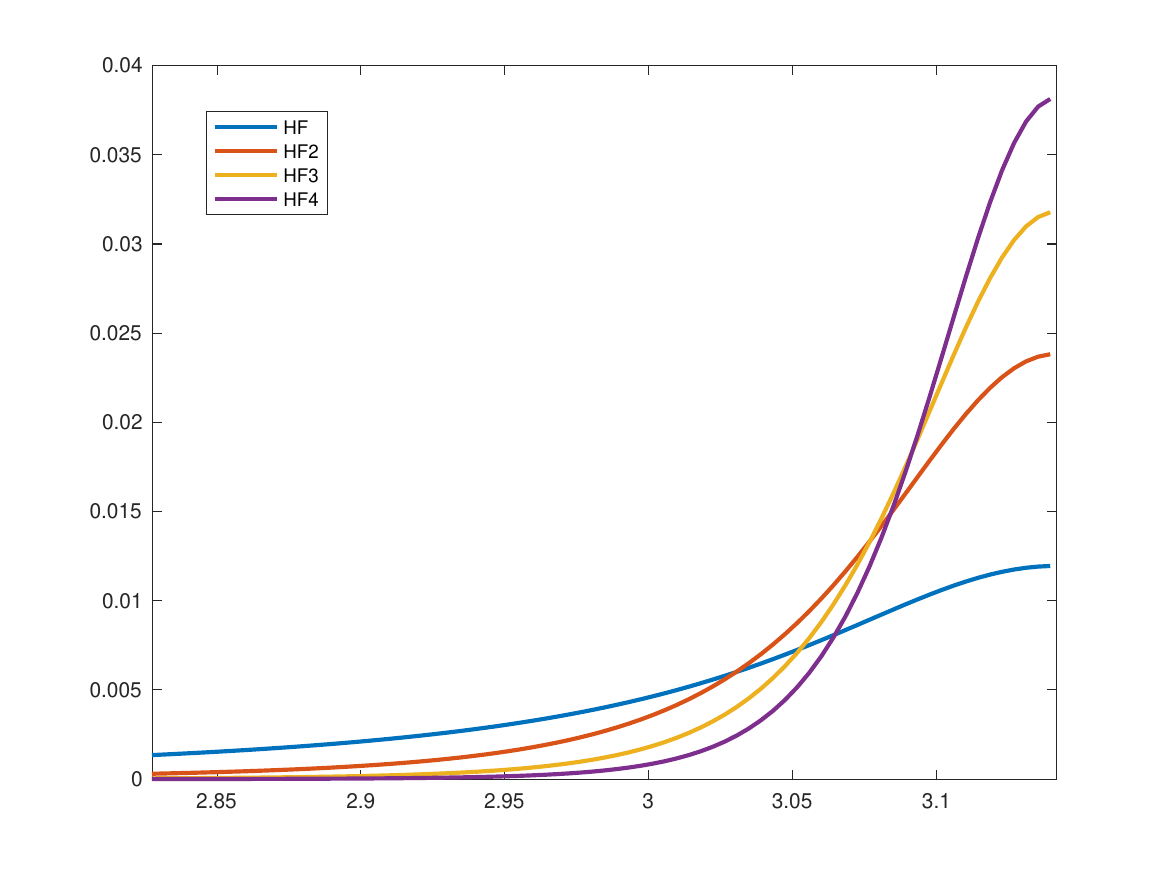}
\caption{Power spectral density of the stationary part of HF and HF$\delta$, with $\delta=2\ldots 4$, with $\beta=0.8$. All those power spectral densities are normalized to one in order to ease the comparison.}
\label{fig_freqHF}
\end{figure}

Finally, also HF$\delta$ and HC$\delta$ are ECLS kernels. Figure \ref{fig_freqHF}
 shows the power spectral density of the stationary part of HF and HF$\delta$ for $\delta=2\ldots 4$.
The higher $\delta$ is, the  more the statistical power is concentrated for  frequencies close to $\pi$.

In order to test the performance of the TC$\delta$ and DC$\delta$ kernels we consider a Monte Carlo study composed by 200 experiments. In each experiment the models and the data are generated likewise to the first Monte Carlo study of Section \ref{sec_GC}, but the input $u$ is a realization drawn from a Gaussian noise with band [0, 0.2]. We consider the following additional estimators for the impulse response:  
\begin{itemize}
\item $\hat g_{DC\delta}$ is the estimator in (\ref{def_RELS}) using the DC$\delta$ kernel (\ref{KDCd}) with $\delta=2\ldots 6$;
\item $\hat g_{TC\delta}$ is the estimator in (\ref{def_RELS}) using the TC$\delta$ kernel (\ref{kernelTCd}) with $\delta=2\ldots 6$.
\end{itemize} 
 \begin{figure*}[ht]
\centering
\includegraphics[width=\textwidth]{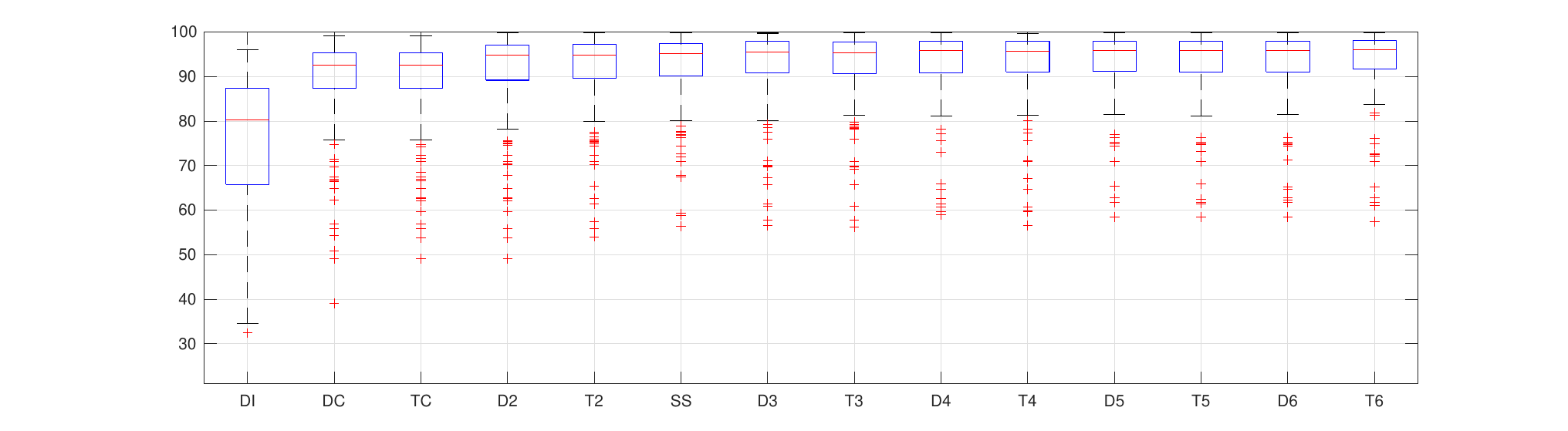}
\caption{Average impulse response fit in the Monte Carlo study composed by 200 experiments.}
\label{fig_sim1}
\end{figure*} 
Figure \ref{fig_sim1} shows the boxplot 
of $\mathrm{AIRF}$ for the estimators. The best one  is $\hat g_{TC6}$, while $\hat g_{DI}$, $\hat g_{DC}$ and  $\hat g_{TC}$ are the worst ones. This result is not surprising because the impulse responses in this Monte Carlo study are enough smooth, see Figure \ref{fig_realsimhl} (top) and, indeed, the TC6 kernel induces more smoothness than the others.

\section{Conclusions}\label{sec_conc}
We have introduced a second-order extension to TC and DC kernels called TC2 and DC2, respectively. 
The latter induces more smoothness than the former.  This idea can be also extended to higher-orders. We also have introduced  a generalized-correlated (GC) kernel which incorporates the DI, DC, TC kernels, i.e. the most popular kernels in system identification, and the DC2 and TC2 kernels. 
We have derived the closed form expression for the determinant and the Cholesky factorization of the inverse matrix of TC2, DC2 and GC. Accordingly, the latter allow to design efficient algorithms for minimizing the negative log-likelihood. In particular, since TC2 and SS kernels produce similar performances for estimating the impulse response, TC2 represents an appealing alternative to SS because it admits an efficient implementation for searching the optimal hyperparameters through marginal likelihood. Finally,
 we have also shown that these new kernels are exponentially convex local stationary and thus it is possible to understand easily their frequency properties.


\begin{thebibliography}{10}

\bibitem{AKAIKE_1974}
H.~Akaike.
\newblock {A new look at the statistical model identification}.
\newblock {\em IEEE Transactions on Automatic Control}, 19(6):716--723, Dec.
  1974.

\bibitem{doi:10.1137/19M1267349}
Martin~S. Andersen and Tianshi Chen.
\newblock Smoothing splines and rank structured matrices: Revisiting the spline
  kernel.
\newblock {\em SIAM Journal on Matrix Analysis and Applications},
  41(2):389--412, 2020.

\bibitem{CARLI2014}
Francesca~P. Carli.
\newblock On the maximum entropy property of the first-order stable spline
  kernel and its implications.
\newblock In {\em IEEE Conference on Control Applications (CCA)}, pages
  409--414, 2014.

\bibitem{carli2012efficient}
Francesca~P Carli, Alessandro Chiuso, and Gianluigi Pillonetto.
\newblock Efficient algorithms for large scale linear system identification
  using stable spline estimators.
\newblock {\em IFAC Proceedings Volumes}, 45(16):119--124, 2012.

\bibitem{carli2011maximum}
Francesca~P Carli, Augusto Ferrante, Michele Pavon, and Giorgio Picci.
\newblock A maximum entropy solution of the covariance extension problem for
  reciprocal processes.
\newblock {\em IEEE Transactions on Automatic Control}, 56(9):1999--2012, 2011.

\bibitem{carli2013efficient}
Francesca~P Carli, Augusto Ferrante, Michele Pavon, and Giorgio Picci.
\newblock An efficient algorithm for maximum entropy extension of
  block-circulant covariance matrices.
\newblock {\em Linear Algebra and its Applications}, 439(8):2309--2329, 2013.

\bibitem{7495008}
Francesca~Paola Carli, Tianshi Chen, and Lennart Ljung.
\newblock Maximum entropy kernels for system identification.
\newblock {\em IEEE Transactions on Automatic Control}, 62(3):1471--1477, 2017.

\bibitem{CHEN2018109}
T.~Chen.
\newblock On kernel design for regularized {LTI} system identification.
\newblock {\em Automatica}, 90:109--122, 2018.

\bibitem{EST_TF_REVISITED_2012}
T.~Chen, H.~Ohlsson, and L.~Ljung.
\newblock On the estimation of transfer functions, regularizations and gaussian
  processes-revisited.
\newblock {\em Automatica}, 48(8):1525--1535, 2012.

\bibitem{chen2018continuous}
Tianshi Chen.
\newblock Continuous-time dc kernel—a stable generalized first-order spline
  kernel.
\newblock {\em IEEE Transactions on Automatic Control}, 63(12):4442--4447,
  2018.

\bibitem{chen2021semiseparable}
Tianshi Chen and Martin~S Andersen.
\newblock On semiseparable kernels and efficient implementation for regularized
  system identification and function estimation.
\newblock {\em Automatica}, 132:109682, 2021.

\bibitem{chen2018regularized}
Tianshi Chen, Martin~S Andersen, Biqiang Mu, Feng Yin, Lennart Ljung, and S~Joe
  Qin.
\newblock Regularized {LTI} system identification with multiple regularization
  matrix.
\newblock {\em Ifac-papersonline}, 51(15):180--185, 2018.

\bibitem{chen2016maximum}
Tianshi Chen, Tohid Ardeshiri, Francesca~P Carli, Alessandro Chiuso, Lennart
  Ljung, and Gianluigi Pillonetto.
\newblock Maximum entropy properties of discrete-time first-order stable spline
  kernel.
\newblock {\em Automatica}, 66:34--38, 2016.

\bibitem{chen2013implementation}
Tianshi Chen and Lennart Ljung.
\newblock Implementation of algorithms for tuning parameters in regularized
  least squares problems in system identification.
\newblock {\em Automatica}, 49(7):2213--2220, 2013.

\bibitem{CHIUSO_PILLONETTO_SPARSE_2012}
A.~Chiuso and G.~Pillonetto.
\newblock A {B}ayesian approach to sparse dynamic network identification.
\newblock {\em Automatica}, 48(8):1553--1565, 2012.

\bibitem{dempster1972covariance}
Arthur~P Dempster.
\newblock Covariance selection.
\newblock {\em Biometrics}, 28:157--175, 1972.

\bibitem{dinuzzo2015kernels}
Francesco Dinuzzo.
\newblock Kernels for linear time invariant system identification.
\newblock {\em SIAM Journal on Control and Optimization}, 53(5):3299--3317,
  2015.

\bibitem{dym1981extensions}
Harry Dym and Israel Gohberg.
\newblock Extensions of band matrices with band inverses.
\newblock {\em Linear algebra and its applications}, 36:1--24, 1981.

\bibitem{matriciToeplitz}
Neville~J. Ford, Dmitry~V. Savostyanov, and Nickolai~L. Zamarashkin.
\newblock On the decay of the elements of inverse {T}riangular toeplitz
  matrices.
\newblock {\em SIAM Journal on Matrix Analysis and Applications},
  35(4):1288--1302, 2014.

\bibitem{FUJIMOTO21}
Y~Fujimoto.
\newblock Efficient implementation of kernel regularization based on {ADMM}.
\newblock In {\em SYSID}, 2021.

\bibitem{gohberg1993classes}
Israel Gohberg, Seymour Goldberg, and Marius~A Kaashoek.
\newblock {\em Classes of linear operators}, volume~63.
\newblock Birkh{\"a}user, 1993.

\bibitem{INFINITE_MATRICES}
P.~Jorgesen, K.~Kornelson, and K.~Shuman.
\newblock {\em {I}terated Function Systems, Moments, and Transformations of
  Infinite Matrices}.
\newblock {A}merican Mathematical {S}ociety, 2011.

\bibitem{LJUNG_SYS_ID_1999}
L.~Ljung.
\newblock {\em System Identification: Theory for the User}.
\newblock Prentice Hall, New Jersey, 1999.

\bibitem{marconato2017filter}
Anna Marconato, Maarten Schoukens, and Johan Schoukens.
\newblock Filter-based regularisation for impulse response modelling.
\newblock {\em IET Control Theory \& Applications}, 11(2):194--204, 2017.

\bibitem{PILLONETTO_DENICOLAO2010}
G.~Pillonetto and G.~De~Nicolao.
\newblock A new kernel-based approach for linear system identification.
\newblock {\em Automatica}, 46:81--93, 2010.

\bibitem{KERNEL_METHODS_2014}
G.~Pillonetto, F.~Dinuzzo, T.~Chen, G.~De~Nicolao, and L.~Ljung.
\newblock Kernel methods in system identification, machine learning and
  function estimation: A survey.
\newblock {\em Automatica}, 50(3):657--682, 2014.

\bibitem{6160606}
Gianluigi Pillonetto and Giuseppe De~Nicolao.
\newblock Kernel selection in linear system identification part i: A gaussian
  process perspective.
\newblock In {\em 50th IEEE Conference on Decision and Control and European
  Control Conference}, pages 4318--4325, 2011.

\bibitem{RASMUSSEN_WILLIAMNS_2006}
C.~Rasmussen and C.~Williams.
\newblock {\em {Gaussian Processes for Machine Learning}}.
\newblock The MIT Press, 2006.

\bibitem{SCHWARZ_1978}
G.~Schwarz.
\newblock {Estimating the Dimension of a Model}.
\newblock {\em The Annals of Statistics}, 6(2):461--464, Mar. 1978.

\bibitem{SODERSTROM_STOICA_1988}
T.~S\"{o}derstr\"{o}m and P.~Stoica.
\newblock {\em System Identification}.
\newblock Prentice-Hall International, Hemel Hempstead, UK, 1989.

\bibitem{wahba1990spline}
Grace Wahba.
\newblock {\em Spline models for observational data}.
\newblock SIAM, 1990.

\bibitem{BSL_CDC}
M.~Zorzi and A.~Chiuso.
\newblock A {B}ayesian approach to sparse plus low rank network identification.
\newblock In {\em Proceedings of the IEEE Conference on Decision and Control},
  pages 7386--7391, Osaka, 2015.

\bibitem{BSL}
M.~Zorzi and A.~Chiuso.
\newblock Sparse plus low rank network identification: A nonparametric
  approach.
\newblock {\em Automatica}, 76:355--366, 2017.

\bibitem{zorzi2020new}
Mattia Zorzi.
\newblock A new kernel-based approach for spectral estimation.
\newblock In {\em European Control Conference (ECC)}, pages 534--539, 2020.

\bibitem{ZORZI2018125}
Mattia Zorzi and Alessandro Chiuso.
\newblock The harmonic analysis of kernel functions.
\newblock {\em Automatica}, 94:125--137, 2018.

\end{thebibliography}
\end{document}